\documentclass{amsart}
\usepackage{latexsym,amsxtra,amscd,ifthen}
\usepackage{amsfonts}
\usepackage{verbatim}
\usepackage{amsmath}
\usepackage{amsthm}
\usepackage{amssymb}

\numberwithin{equation}{section}

\newcommand{\mattwo}[1]{\left[\begin{array}{rr} #1 \end{array}\right]}
\newcommand{\matthree}[1]{\left[\begin{array}{rrr} #1 \end{array}\right]}
\newcommand{\matfour}[1]{\left[\begin{array}{rrrr} #1 \end{array}\right]}
\newcommand{\matfive}[1]{\left[\begin{array}{rrrrr} #1 \end{array}\right]}

\newcommand{\mateight}[1]{\left[\begin{array}{rrrrrrrr} #1 \end{array}\right]}
\theoremstyle{plain}
\newtheorem{theorem}{Theorem}[section]
\newtheorem{lemma}[theorem]{Lemma}
\newtheorem{proposition}[theorem]{Proposition}
\newtheorem{corollary}[theorem]{Corollary}
\newtheorem{conjecture}[theorem]{Conjecture}

\theoremstyle{definition}
\newtheorem{definition}[theorem]{Definition}
\newtheorem{example}[theorem]{Example}

\newtheorem{remark}[theorem]{Remark}
\newtheorem{question}[theorem]{Question}

\makeatletter              
\let\c@equation\c@theorem  
\makeatother

\DeclareMathOperator{\hdet}{hdet}

\DeclareMathOperator{\Ext}{Ext}

 \DeclareMathOperator{\tr}{tr}

\DeclareMathOperator{\Aut}{Aut_{gr}}

\DeclareMathOperator{\GKdim}{GKdim}

\newcommand{\pp}{\textup{p}}

\begin{document}

\title[Shephard-Todd-Chevalley Theorem]
{Shephard-Todd-Chevalley Theorem \\
for Skew Polynomial Rings}
\author{E. Kirkman, J. Kuzmanovich and J.J. Zhang}

\address{Kirkman: Department of Mathematics,
P. O. Box 7388, Wake Forest University, Winston-Salem, NC 27109}

\email{kirkman@wfu.edu}

\address{Kuzmanovich: Department of Mathematics,
P. O. Box 7388, Wake Forest University, Winston-Salem, NC 27109}

\email{kuz@wfu.edu}

\address{zhang: Department of Mathematics, Box 354350,
University of Washington, Seattle, Washington 98195, USA}

\email{zhang@math.washington.edu}

\begin{abstract}
We prove the following generalization of the classical
Shephard-Todd-Chevalley Theorem. Let $G$ be a finite group of
graded algebra automorphisms of a skew polynomial ring
$A:=k_{p_{ij}}[x_1,\cdots,x_n]$. Then the fixed subring $A^G$ has
finite global dimension if and only if $G$ is generated by
quasi-reflections. In this case the fixed subring $A^G$ is
isomorphic a skew polynomial ring with possibly different
$p_{ij}$'s. A version of the theorem is proved also for abelian
groups acting on general quantum polynomial rings.
\end{abstract}








\maketitle


\setcounter{section}{-1}
\section{Introduction}
\label{sec0}

The classical Shephard-Todd-Chevalley Theorem states that if $G$ is
a finite group acting faithfully on a finite dimensional
$k$-vector space $V:=\bigoplus_{i=1}^n kx_i$, then the
fixed subring $k[x_1,\cdots,x_n]^G$ is isomorphic to
$k[x_1,\cdots,x_n]$ if and only if $G$ is generated by
pseudo-reflections of $V$. (In the rest
of this paper a pseudo-reflection will be called simply a
reflection). It is natural to ask if a
version of the Shephard-Todd-Chevalley Theorem holds for regular algebras,
which are algebras of finite global dimension that are noncommutative
analogs of commutative polynomial rings. In \cite{KKZ1} the
authors began our investigation of conditions on a group $G$ of graded
algebra automorphisms of a regular algebra $A$ so that $A^G$ is
isomorphic to $A$, or that $A^G$ itself has finite global dimension and hence is regular.

The main goal of this paper is to provide  evidence that
a noncommutative version of Shephard-Todd-Chevalley Theorem should
hold.  Let $k$ be a base commutative
field of characteristic zero, and let $n$ be a positive integer.
Let $\{p_{ij}\;|\; 1\leq i<j\leq n\}$ be a set of nonzero
scalars in $k$. The skew polynomial ring $A:=k_{p_{ij}}[x_1,\cdots,
x_n]$ is defined to be the $k$-algebra generated by $x_1,\cdots,
x_n$ and subject to the relations
$$x_jx_i=p_{ij}x_ix_j$$
for all $1\leq i<j\leq n$. The algebra $A$ is a connected
${\mathbb N}$-graded algebra with $\deg x_i=1$ for all $i$. Let
$\Aut(A)$ be the group of all graded algebra automorphisms
of $A$.  Our main result is the following theorem.

\begin{theorem}[Theorem \ref{xxthm4.5}]
\label{xxthm0.1}
Let $A=k_{p_{ij}}[x_1,\cdots, x_n]$ and
$G$ a finite subgroup of $\Aut(A)$. Then the fixed subring
$A^G$ has finite global dimension if and only if $G$ is
generated by quasi-reflections.
\end{theorem}

When the fixed subring $A^G$ has finite global dimension, it need not be isomorphic to
$A$ itself (as in the classical Shephard-Todd-Chevalley Theorem), but it is isomorphic to another skew polynomial
ring $k_{p'_{ij}}[x_1,\cdots,x_n]$. The definition of a
quasi-reflection is given in Definition \ref{xxdefn1.2}.
As we will see in Section 1  most  quasi-reflections
are reflections in the classical sense,
with only one interesting exception, called mystic
reflections [Definition \ref{xxdefn1.3}(b)].

If $p_{ij}\neq \pm 1$ for all $i<j$, then the quasi-reflections
of $A$ are easy to describe and the proof of
Theorem \ref{xxthm0.1} is simple. The complication
appears only when some of $p_{ij}$s are $\pm 1$.
The idea in the proof is to group together certain
quasi-reflections so that we have a partition on the space
$\bigoplus_{i=1}^n kx_i$ under certain rules determined by
the parameter set $\{p_{ij}\}$ and the group $G$. Using
this partition we can reduce the question to each block, and
within each block the algebra is either commutative or
connected by mystic reflections. We believe that this kind
of partition should be useful for proving a more general
version of the Shephard-Todd-Chevalley Theorem conjectured
as follows.

\begin{conjecture}
\label{xxcon0.2}
Let $B$ be a quantum polynomial ring
[Definition \ref{xxdefn1.1}(2)] and let $G$ be a finite subgroup
of $\Aut(B)$. Then  $B^G$ has finite global dimension
if and only if $G$ is generated by quasi-reflections.
\end{conjecture}

Other evidence for the above conjecture is the following
theorem, which settles the case when $G$ is abelian.

\begin{theorem}[Theorem \ref{xxthm5.3}]
\label{xxthm0.3} Let $B$ be a quantum polynomial ring
and let $G$ be a finite abelian subgroup of $\Aut(B)$.
Then $B^G$ has finite global dimension if and only if
$G$ is generated by quasi-reflections.
\end{theorem}

Theorem \ref{xxthm0.3} is a generalization of
\cite[Theorem 0.5]{KKZ1}. The combination of
Theorem \ref{xxthm0.1} with Theorem \ref{xxthm0.3}
covers only a small part of Conjecture
\ref{xxcon0.2}. For a different reason (see
Proposition \ref{xxprop0.5}) we also verify the
Shephard-Todd-Chevalley Theorem for the quantum
$2\times 2$ matrix algebra ${\mathcal O}_q(M_2)$,
which has some non-abelian reflection groups
[Proposition \ref{xxprop5.8}].

The new quasi-reflections possessed by some regular algebras
suggest extending the notion of ``reflection group" to include
groups $G$ with a representation by automorphisms of
a regular algebra $A$ giving a regular fixed subring $A^G$.
This leads to a few natural questions:

\begin{question}
\label{xxque0.4}
\begin{enumerate}
\item
Can every ``reflection group'' of a noncommutative regular
algebra be realized as a reflection group in the classical sense?
\item
Is there a version of the Shephard-Todd-Chevalley Theorem
for finite quantum group (i.e., finite dimensional
Hopf algebra) actions?
\item
Using Hopf algebra actions, instead of group actions,
can we obtain more regular fixed subrings of a given quantum
polynomial ring?
\end{enumerate}
\end{question}

A secondary goal of this paper is to try to answer the above
questions. J. Alev asked a question similar to
Question \ref{xxque0.4}(c). In Example \ref{xxex6.2} we show
that some groups of mystic reflections are not reflection groups
in the classical sense, which answers Question \ref{xxque0.4}(a)
negatively. In Section 5, we prove the following result.

\begin{proposition}
\label{xxprop0.5} [Propositions \ref{xxprop5.7} and
\ref{xxprop5.8}]
Let $C$ be the quantum $2\times 2$ matrix algebra
${\mathcal O}_q(M_2)$. Suppose $q\neq \pm 1$.
\begin{enumerate}
\item[(i)]
For every finite group $G\subset \Aut(C)$,
$C^G$ is not isomorphic to a skew polynomial ring
(as defined before Theorem \ref{xxthm0.1}).
For every non-trivial finite group $G\subset \Aut(C)$,
$C^G$ is not isomorphic to $C$.
\item[(ii)]
Suppose $q$ is a root of $1$. Then there is a finite dimensional
Hopf algebra $H$ with a Hopf action on $C$ such that $C^H$ is
isomorphic to a skew polynomial ring. As a consequence, $C^H$ has
finite global dimension and is not isomorphic to $C^G$
for any finite group $G\subset \Aut(C)$.
\end{enumerate}
\end{proposition}

Hence by Proposition \ref{xxprop0.5}(ii) we do obtain more regular
fixed subrings using Hopf
algebra actions, which answers Question \ref{xxque0.4}(c).
As far as Question \ref{xxque0.4}(b), Proposition \ref{xxprop5.7}
suggests that there should be some version of
Shephard-Todd-Chevalley Theorem for Hopf algebra actions,
however, we do not have an explicit conjectural statement.
Some related work about Hopf algebra actions on Artin-Schelter
Gorenstein algebras will be presented in \cite{KKZ2}.

The paper is organized as follows: Some preliminary material
is reviewed in Section 1. Sections 2 and 3 contain some
analysis of quasi-reflections of skew polynomial rings.
Theorem \ref{xxthm0.1} is proved in Section 4 and
Theorem \ref{xxthm0.3} is proved in Section 5. Section 6
contains some remarks about mystic reflection groups.

\section{Definitions}
\label{sec1}

Throughout $k$ is a base field of characteristic zero.
In \cite{KKZ1} the authors assume that $k$ is algebraically
closed. For simplicity and our convenience we continue to
assume that $k$ is algebraically closed since we will use
several results from \cite{KKZ1}, though
all the main assertions in this paper hold without that assumption.
The opposite ring of an algebra $A$ is denoted by $A^{op}$.
Usually we are working with
left $A$-modules, and a right $A$-module can be viewed as a
left $A^{op}$-module.

An algebra $A$ is called {\it connected ${\mathbb N}$-graded} or
{\it connected graded} if
$$A=k\oplus A_1\oplus A_2\oplus \cdots$$
and $A_iA_j\subset A_{i+j}$ for all $i,j\in {\mathbb N}$.
The Hilbert series of $A$ is defined to be
$$H_A(t)=\sum_{i\in {\mathbb Z}} (\dim A_i)t^i.$$
The class of algebras considered in this paper are the graded algebras with finite global dimension.
Sometimes we need the following more restrictive class of quantum polynomial rings.

\begin{definition}
\label{xxdefn1.1}
Let $A$ be a connected graded algebra.

(1) We call $A$ {\it Artin-Schelter Gorenstein}
if the following conditions hold:
\begin{enumerate}
\item
$A$ has graded injective dimension $d<\infty$ on
the left and on the right,
\item
$\Ext^i_A(k,A)=\Ext^i_{A^{op}}(k,A)=0$ for all
$i\neq d$, and
\item
$\Ext^d_A(k,A)\cong \Ext^d_{A^{op}}(k,A)\cong k(l)$
for some
$l$.
\end{enumerate}
If in addition,
\begin{enumerate}
\item[(d)]
$A$ has finite (graded) global dimension, and
\item[(e)]
$A$ has finite Gelfand-Kirillov dimension,
\end{enumerate}
then $A$ is called {\it Artin-Schelter regular}
(or {\it regular} for short) of dimension $d$.

(2) If $A$ is a noetherian, regular graded domain of global
dimension $n$ and $H_A(t)=(1-t)^{-n}$,
then we call $A$ {\it a quantum polynomial ring}
of dimension $n$.
\end{definition}

Skew polynomial rings (with $\deg x_i=1$) are quantum polynomial
rings. Next we recall some definitions of the noncommutative
versions of reflections from \cite{KKZ1}. Let $g\in \Aut(A)$,
the graded algebra automorphisms of $A$, then
the trace function of $g$ is defined to be
$$Tr_A(g,t)=\sum_{i=0}^{\infty} \tr(g|_{A_i}) t^i
\in k[[t]],$$
where $\tr(g|_{A_i})$ is the trace of the linear map
$g|_{A_i}$.
Note that $Tr_A(g,0)=1$.

\begin{definition}
\label{xxdefn1.2}
\cite[Definition 2.2]{KKZ1}
Let $A$ be a regular graded algebra such that
$$H_A(t) =\frac{1}{(1-t)^{n}f(t)}$$
where $f(1)\neq 0$ (so that $n=\GKdim A$).
Let $g$ be a graded algebra automorphism of $A$.
We say that $g$ is a {\it quasi-reflection} of $A$ if
$$Tr_A(g,t)={\frac{1}{(1-t)^{n-1}q(t)}}$$
for $q(1) \neq 0$. If $A$ is a quantum polynomial ring,
then $H_A(t)=(1-t)^{-n}$. In this case $g$ is a
{\it quasi-reflection} if and only if
$$Tr_A(g,t)={\frac{1}{(1-t)^{n-1}(1-\lambda t)}}$$
for some scalar $\lambda \neq 1$. (Note that we have chosen not
to call the identity map a quasi-reflection).
\end{definition}

Quasi-reflections of quantum polynomial rings are
characterized in \cite[Theorem 3.1]{KKZ1}.

\begin{definition}
\label{xxdefn1.3} Let $A$ be a quantum polynomial ring of
global dimension $n$. If $g\in \Aut(A)$ is a quasi-reflection
of finite order, by \cite[Theorem 3.1]{KKZ1} $g$ is described in one of
the following two cases.
\begin{enumerate}
\item
There is a basis of $A_1$, say $\{y_1,\cdots,y_n\}$,
such that $g(y_j)=y_j$ for all $j\leq n-1 $ and $g(y_n)
=\lambda y_n$. Namely, $g|_{A_1}$ is a reflection.
In this case $g$ is called a {\it reflection} of $A$,
$y_n$ is called a {\it non-invariant eigenvector}
 of $g$,
and $E_g:=\oplus_{i<n} ky_i$ is called the {\it invariant
eigenspace} of $g$.
\item
The order of $g$ is $4$, and there is a basis of $A_1$,
say $\{y_1,\cdots,y_n\}$,
such that $g(y_j)=y_j$ for all $j\leq n-2$, $g(y_{n-1})
=i\; y_{n-1}$, and $g(y_n)=-i\; y_n$ (where $i^2=-1$).
In this case $g$ is called a {\it mystic reflection}
of $A$, $\{y_{n-1},y_n\}$ is called a pair of {\it non-invariant
eigenvectors} of $g$, and $E_g:=\oplus_{i<n-1} ky_i$ is
called the {\it invariant eigenspace} of $g$.
\end{enumerate}
\end{definition}

The homological determinant $\hdet$ of $g$ is defined in
\cite[Definition 2.3]{JoZ}. We refer to \cite{JoZ} for some
details. We use $\det$ for the usual determinant of a $k$-linear map.

\begin{proposition}
\label{xxprop1.4} Let $A$ be a quantum polynomial ring and
let $g$ be a quasi-reflection of $A$.
\begin{enumerate}
\item
If $g$ is a reflection, then the non-invariant eigenvector $y_n$ is
a normal element of $A$, and $\det g|_{A_1}=\hdet g=\lambda
\neq 1$.
\item
Suppose $g$ is a mystic reflection with a pair of
non-invariant eigenvectors $y_{n-1}$ and $y_n$.  Then
$E_g y_{n-1}=y_{n-1}E_g$, $E_g y_{n}=y_{n}E_g$,
$\det g|_{A_1}=1$, and $\hdet g=-1$.
\item
If $g$ is a mystic reflection, then $y_{n-1}^2=y_n^2$ up to a scalar,
and there are two elements
$y'_{n-1},y'_n$ with $ky'_{n-1} + k y'_{n} =  ky_{n-1}+ky_n$ and $y'_ny'_{n-1}=-y'_{n-1}y'_n$.
Furthermore, $g(y'_{n-1})=y'_n$ and $g(y'_n)=-y'_{n-1}$.
\item
Let $G$ be generated by mystic reflections. Then $G$ does not
contain any reflections.
\end{enumerate}
\end{proposition}

\begin{proof} (a) The first assertion follows from the proof
of \cite[Lemma 5.2(a,d)]{KKZ1}. The second assertion is clear
from the definition.

(b) The first assertion follows from the proof of
\cite[Proposition 4.3]{KKZ1}.

It is clear that if $g$ is a mystic reflection then
$\det g|_{A_1}=1$. By \cite[Lemma 4.2]{KKZ1},
$\hdet g=-1$.

(c) The first two assertions are  \cite[Proposition 4.3(c)]{KKZ1}.
For the last assertion we set $y'_{n-1}=y_{n-1}+i y_n$
and $y'_n=g(y'_{n-1})=iy_{n-1}+y_n$.

(d) By part (b), $\det f|_{A_1}=1$ for all $f\in G$. By
part (a), such an $f$ cannot be a reflection.
\end{proof}

\begin{proposition}
\label{xxprop1.5}
Let $A$ be a noetherian regular algebra and
$G$ be a finite group of automorphisms of $A$.
\begin{enumerate}
\item
Let $R$ be the subgroup of $G$ generated by
quasi-reflections. Then $R$ is a normal subgroup of $G$,
and if both $A^R$ and $A^G$ have finite
global dimension, then $R=G$.
\item
Suppose, for every finite $H\subset \Aut(A)$ generated by
quasi-reflections, $A^H$ has finite global dimension.
Then $A^G$ has  finite global dimension if and only if
$G$ is generated by quasi-reflections.
\end{enumerate}
\end{proposition}

\begin{proof} (a) By \cite[Lemma 1.10(c)]{KKZ1}, if $A^G$
has finite global dimension, then it is regular.

Let $n=\GKdim A$. Since the property of $h$ being a
quasi-reflection is determined by $Tr(h,t)$ and
$Tr(g h g^{-1},t)=Tr(h,t)$ for any $g\in G$, then $R$ is a
normal subgroup of $G$. Thus the quotient group $G/R$
acts on the fixed subring $A^R$, and the fixed subring
$(A^R)^{G/R}$ equals $A^G$.
Let $g \in G \setminus R$. Then by \cite[Lemma 5.2]{JiZ}
$$ Tr_{A^R}(g|_{A^R},t) = \frac{1}{|R|} \sum_{h \in R}
Tr_A(hg,t).$$
Note that no $hg$ can be a quasi-reflection, for otherwise
we would have $g \in R$. Hence $Tr_A(hg,t)$ must have the form
$$ Tr_A(hg,t) = \frac{1}{(1-t)^m p(t)}$$
where $m \leq n-2$ and $p(1) \neq 0$.  This means that
$(1-t)^{n-2} Tr_{A^R}(g|_{A^R},t)$ is analytic at $t=1$ and
$g|_{A^R}$ cannot be a quasi-reflection. Since $(A^R)^{G/R}$
is regular, $G/R$ must contain a quasi-reflection by
\cite[Theorem 2.4]{KKZ1}.  Hence $G/R$ is trivial and $G=R$.

(b) This is immediate from part (a).
\end{proof}

It follows from this proposition that to prove a result like Theorem \ref{xxthm0.1}
we need only to show one direction, namely, to show that if $G$
is generated by quasi-reflections then $A^G$ has finite global
dimension.

In the rest of the section we review two elementary lemmas.

\begin{lemma}
\label{xxlem1.6}
Let $A(i)$, for $i=1,\cdots,d$, be regular algebras
(or connected graded algebras having finite global dimension).
Then the tensor product $A(1)\otimes_k \cdots \otimes_k A(d)$
is a regular algebra (or an algebra having finite global
dimension).
\end{lemma}

\begin{proof} By induction we may assume $d=2$.
Then we use the proof of \cite[Proposition 4.5(a)]{YZ}.
\end{proof}

The algebra $A(1)\otimes_k \cdots \otimes_k A(d)$ in the above
lemma is connected ${\mathbb N}^d$-graded. Let $B$ be any
${\mathbb N}^d$-graded algebra. Let $\phi:=
\{\phi_s\;|\; s=1,\cdots,d\}$
be a sequence of commuting ${\mathbb N}^d$-graded algebra
automorphisms of $B$. Then we can define a twisted algebra,
as in \cite{Z}, $B^\phi$ as follows: $B^\phi=B$
as a ${\mathbb N}^d$-graded vector space, and the new
multiplication $*$ of $B^\phi$ is determined by
$$ a* b= a \phi^{|a|}(b)$$
where $\phi^{|a|}=\phi_1^{a_1}\cdots \phi_d^{a_d}$ if the
degree of $a$ is $|a|=(a_1,\cdots,a_d)$. We refer to
\cite{Z} for some basic properties of twisted algebras.
The following lemma is elementary.

\begin{lemma}
\label{xxlem1.7}
Let $B$ be a connected ${\mathbb N}^d$-graded
algebra that is regular as a connected graded algebra.
Let $\phi:=\{\phi_s\;|\; s=1,\cdots,d\}$
be a sequence of commuting ${\mathbb N}^d$-graded algebra
automorphisms of $B$.
\begin{enumerate}
\item
The algebra $B^\phi$ is a connected ${\mathbb N}^d$-graded
algebra that is regular as a connected graded algebra.
\item
Suppose $B$ is a skew polynomial ring with generators $\{x_i\}$.
($B$ may not be generated in degree 1.)
If  each $x_i$ is ${\mathbb N}^d$-homogeneous and if $\phi_s$ maps
$x_i$ to $q_{si} x_i$ for all $s$ and $i$, where $q_{si}$
are some nonzero scalars in $k$, then $B^\phi$ is
isomorphic to a skew polynomial ring.
\end{enumerate}
\end{lemma}

\section{Block decomposition, elementary transformations\\
and quasi-reflections of skew polynomial rings}
\label{sec2}

In Sections 2, 3 and 4 we fix an integer $n\geq 2$ and
a set of nonzero scalars $\{p_{ij}\;|\; 1\leq i<j\leq n\}$.
Let $A$ be the skew polynomial ring $k_{p_{ij}}[x_1,\cdots,x_n]$
that is generated by $\{x_1,\cdots,x_n\}$ and subject to
the relations
$x_jx_i=p_{ij}x_ix_j$
for all $1\leq i<j\leq n$. We will try to use other letters
for a general algebra. For $i\geq j$ set
\begin{equation}
\label{2.0.1}
p_{ij}=\begin{cases}
1& i=j\\ p_{ji}^{-1} &i>j.\end{cases}
\tag{2.0.1}
\end{equation}
Let $[n]$ be the set $\{1,2,\cdots, n\}$ and let $\pp$ be
the set $\{p_{ij}\;|\; i,j\in [n]\}$. Following \eqref{2.0.1}
we have $x_jx_i=p_{ij}x_ix_j$ for all $i,j\in [n]$.
Let $A_1=\oplus_{i=1}^{n} k x_i$ which is called the
{\it generating space} of $A$. Let $\{x'_1,\cdots,x'_n\}$
be another basis of $A_1$. We call $\{x'_1,\cdots,x'_n\}$
a $\pp$-basis if $x'_jx'_i=p_{ij}x'_ix'_j$ holds in $A$
for all $i,j\in [n]$.

\begin{definition}
\label{xxdefn2.1}
Fix a set $\pp$ satisfying \eqref{2.0.1} and fix an $i\in [n]$.
\begin{enumerate}
\item
We define the {\it block} containing $i$ to be
$$B(i)=\{i'\in [n]\;|\; p_{ij}=p_{i'j}\quad \forall j\in [n]\}.$$
We say that $i$ and $j$ are in the same block if $B(i)=B(j)$.
\item
We use the blocks to define an equivalence relation
on $[n],$ and then $[n]$ is a disjoint union of blocks
$$[n]=\bigcup_{i\in W} B(i)$$
for some index set $W$. The above equation is called the {\it block
decomposition}.
\item
Let $m:=|W|$, which is no more than $n$. We call $m$ the {\it number}
of blocks. For simplicity, we may identity $W$ with $[m]$. Also
we will write $B(i)$ as $B_w$ for some $w\in [m]$.
\end{enumerate}
\end{definition}

The block decomposition is a partition of $[n]$ and it is
uniquely determined by $\pp$. By abusing the language, we also
call the subspace $\bigoplus_{s\in B(i)} kx_s$ a block.

In the lemma below we will rescale using a square root.  In
this and in other times that we do this rescaling, we choose a particular
square root $q_{ij}=\sqrt{p_{ij}}$ for $i < j$, and then take
$q_{ji} = 1/q_{ij}$, using the square root chosen for $q_{ij}$.
\begin{lemma}
\label{xxlem2.2}
Let $W=[m]$ and let $\{B_w\;|\; w\in [m]\}$ be the set of blocks.
Let $P_w=\oplus_{i\in B_w} kx_i$ for each $w$.
\begin{enumerate}
\item
If $i,j\in B_w$, then $p_{ij}=1$. This means that the subalgebra
generated by $x_i$ for all $i\in B_w$ is commutative.
As a consequence, if
$p_{ij}\neq 1$ for all $i<j$, then $W=[n]$ and $B(i)=\{i\}$
for all $i\in [n]$.
\item
The algebra $A$ is an ${\mathbb N}^m$-graded algebra
with $\deg x_i=(0,\cdots, 0,1,0,\cdots,0)$ where
$1$ is in the $w$th position if $i\in B_w$.
\item
For every $w\in W$, pick any $k$-linear basis
$\{x'_i\;|\; i\in B_w\}$ for $P_w$,
then there is an ${\mathbb N}^m$-graded algebra automorphism
of $A$ determined by the map
$\theta: x_i\to x'_i$ for all $i\in [n]$.
\item
For each $w\in W$, pick $j\in B_w$.
Define the automorphism $\phi_w:
x_i\to q_{ji} x_i$ for all $i$. (Recall
that $q_{ij}=\sqrt{p_{ij}}$.)
Then $\phi=\{\phi_w\}$ is a sequence of commuting
${\mathbb N}^m$-graded algebra automorphisms of $A$.
\item
The graded twist  $A^{\phi}$ is the commutative
polynomial ring $k[x_1,\cdots,x_n]$.
\end{enumerate}
\end{lemma}

\begin{proof} (a,b,c,d) These are straightforward.

(e) The new multiplication of the graded twist $A^\phi$ is
determined by
$$x_s* x_t=x_s \phi_w(x_t)=x_s q_{st} x_t$$
for $s\in B_w$.
Then
$$x_t*x_s=x_t q_{ts} x_s=p_{st}x_sx_t q_{ts}
= q_{st} x_s x_t=x_s* x_t.$$
Therefore $A^\phi$ is commutative.
\end{proof}

We call the algebra automorphism in Lemma \ref{xxlem2.2}(c)
an {\it elementary transformation}.
There is a class of obvious isomorphisms between two
skew polynomial rings. If there is a permutation
$\sigma\in S_n$ such that $p'_{ij}=p_{\sigma(i)\sigma(j)}$
for all $i,j\in [n]$, then $k_{p'_{ij}}[x_1,\cdots,x_n]$
is isomorphic to $k_{p_{ij}}[x_1,\cdots,x_n]$ by
sending $x_i\to x_{\sigma(i)}$ for all $i\in [n]$.
Such an isomorphism is called a {\it permutation}.

Next we define some standard quasi-reflections of skew
polynomial rings.

\begin{definition}
\label{xxdefn2.3}
As usual we fix a set of scalars
$\pp$ and let $A$ be the skew polynomial ring
$k_{\pp}[x_1,\cdots,x_n]$. Let $\lambda$ be a nonzero
scalar.
\begin{enumerate}
\item
Let $s\in [n]$ and suppose that $\lambda$ is not $1$.
Let $\theta_{s,\lambda}$ be the automorphism of $A$ determined by
$$\theta_{s,\lambda}(x_i)=\begin{cases}
x_i & i\neq s\\ \lambda x_s & i=s.\end{cases}$$
This map is called a {\it standard reflection} of $A$.
\item
Let $s,t\in [n]$. Suppose that $p_{st}=-1$ and that $p_{sj}=p_{tj}$
for all $j\in [n]\setminus \{s,t\}$.
Let $\tau_{s,t,\lambda}$ be the automorphism of $A$ determined by
$$\tau_{s,t,\lambda}(x_i)=\begin{cases}
x_i & i\neq s,t \\ \lambda x_t & i=s\\
-\lambda^{-1} x_s & i=t.\end{cases}$$
This map is called a {\it standard mystic reflection} of $A$.
\item
If $g$ is an elementary transformation, then $g
\theta_{s,\lambda} g^{-1}$ is called an {\it elementary
reflection}.
\item
 If $g$ is an elementary transformation, then $g
\tau_{s,t,\lambda} g^{-1}$ is called an {\it elementary
mystic reflection}.
\end{enumerate}
\end{definition}

The following easy lemma justifies the names given above.

\begin{lemma}
\label{xxlem2.4}
Using the notation above:
\begin{enumerate}
\item
Each $\theta_{s,\lambda}$ is an algebra automorphism of $A$ that is
a reflection of $A$ in the sense of Definition \ref{xxdefn1.3}(a).
\item
Each $\tau_{s,t,\lambda}$ is an algebra automorphism of $A$ that is
a mystic reflection of $A$ in the sense of Definition
\ref{xxdefn1.3}(b).
\item
If $g$ is in $\Aut(A)$ and $\theta$ is a reflection,
then $g\theta g^{-1}$ is a reflection in the sense of
Definition \ref{xxdefn1.3}(a).
\item
If $g$ is in $\Aut(A)$ and $\tau$ is a mystic reflection,
then $g\tau g^{-1}$ is a mystic reflection in the sense of
Definition \ref{xxdefn1.3}(b).
\end{enumerate}
\end{lemma}

In the rest of this section we study the reflections of $A$.
For any $f\in A_1$, we write $f=\sum_{i} a_i x_i$.
If $f\neq 0$, we define $I_f=\{i\;|\; a_i\neq 0\}$.
We also write $f=\sum_{i\in I_f} a_i x_i$.

\begin{lemma}
\label{xxlem2.5}
Let $f$ and $g$ be nonzero elements in $A_1$.
\begin{enumerate}
\item
Suppose $f$ is a normal element of $A$. Then
$p_{ij}=p_{i'j}$ for all $i,i'\in I_f$ and for all $j\in [n]$,
or equivalently, $I_f\subset B_w$ for some $w\in W$.
As a consequence, $p_{ii'}=1$ for all $i,i'\in I_f$.
\item
Suppose $fg=qgf$ for some $1\neq q\in k^{\times}:=k\setminus \{0\}$.
Then $p_{ij}=q$ for all $j\in I_f$ and $i\in I_g$.
As a consequence, $I_f\cap I_g=\emptyset$.
\item
Suppose $fg=gf$. Then $p_{ij}=1$ for all
$j\in I_f$ and $i\in I_g$ and $\{i,j\}\not\subset I_f\cap
I_g$.
\item
Suppose that $p_{ij}\neq 1$ for all $i<j$.
Then every normal element in $A$ is of the form $cx_i$ for
some $c\in k$ and for some $i\in [n]$.
\item
Suppose that $p_{ij}\neq 1$ for all $i<j$.
Let $\phi$ be an automorphism of $A$. Then
there is a permutation $\sigma\in S_n$ and
$c_s\in k^{\times}$ for all $s\in [n]$ such that
$$\phi(x_s)=c_s x_{\sigma(s)}$$ for all $s\in [n]$.
\end{enumerate}
\end{lemma}

\begin{proof} (d) is a consequence of (a), and (e) is
a consequence of (d). So we prove (a,b,c) next.

All three statements are easy to check when
$|I_f|=1$, or $p_{ij}=1$ for all $i$ and $j$,
or $n\leq 2$.

We now assume that (i) $n\geq 3$, (ii) $|I_f|\geq 2$, and
(iii) $p_{ij}\neq 1$ for some $i,j$. Let $A'=A/(x_{i_0})$
for a choice of $i_0\in [n]$, and let $f'$ and $g'$ be
the images of $f$ and $g$ in $A'$, respectively.

(a) Since $|I_f|\geq 2$, $f'$ is nonzero for any $i_0$.
Since $f$ is normal, we have
$$x_{i_0}f=f(\sum_i b_i x_i)$$
which implies
$$0=f'(\sum_{i\neq i_0} b_i x_i)$$
in $A'$. Hence $b_i=0$ for all $i\neq i_0$ and
$$x_{i_0}f=b_{i_0} f x_{i_0}.$$
Therefore $p_{ii_0}=p_{i'i_0}=b_{i_0}$ for all $i,i'\in I_f$.
Since $i_0$ is arbitrary, part (a) is proved.

(b) By symmetry, we may also assume that $|I_g|\geq 2$.
Let $i_0=n$. Write $f=f'+a x_n$ and $g=g'+b x_n$. Then
$f'g'=qg'f'$ in $A'$. If $a=b=0$, then the assertion
follows trivially from $f'g'=qg'f'$ and induction. Now
assume that $a\neq 0$, and so we can assume that
$a=1$ without loss of generality. Note that the
equations $fg=qgf$ and $f'g'=qg'f'$
imply that
$$x_n g'+bf'x_n+ b x_n^2=
q(g'x_n+bx_n f'+bx_n^2).$$
Since $q\neq 1$, then $b=0$ and the assertion follows from
that fact $x_ng'=q g'x_n$.

(c) For any $i\in I_g$ and $j\in I_f$ such that $\{i,j\}
\not\subset I_f\cap I_g$. Since $n\geq 3$, then there is an $i_0$
that is neither $i$ nor $j$. Then we have $f'g'=g'f'$ in
$A'=A/(x_{i_0})$. The assertion follows by induction.
\end{proof}

The following proposition is one of the main results of this section.

\begin{proposition}
\label{xxprop2.6}
Let $\theta$ be a reflection of $A$. Then
$\theta$ is an elementary reflection, namely, up to an
elementary transformation, $\theta=\theta_{i,\lambda}$
for some $i\in [n]$ and for some $1\neq \lambda\in k^{\times}$.
\end{proposition}

\begin{proof} Let $v$ be a non-invariant eigenvector
of $\theta$. Then $\theta(v)=\lambda v$ for some
nonzero scalar $\lambda\neq 1$, and $\theta$ becomes
the identity on $A/(v)$.

By Lemma \ref{xxlem2.5}(a) $I_v\subset B_w$ for some $w\in W$.
We claim that
\begin{enumerate}
\item[(i)]
$\theta(x_s)=x_s$ for all $s\not\in B_w$.
\item[(ii)]
$\theta(x_s)\in \sum_{i\in B_w} kx_i$ for
all $s\in B_w$.
\end{enumerate}
Since $\theta$ becomes the identity on $A/(v)$,
for every $i$,  $\theta(x_i)=x_i+c_i v$ for some
$c_i\in k$. If $s\in B_w$, then $\theta(x_s)
=x_s+c_s v\in \sum_{i\in B_w} kx_i$, establishing condition
(ii). For any $s\not\in B_w$, we claim that $\theta(x_s)=x_s$.
If not, then $\theta(x_s)=x_s+c v$ for some $c\neq 0$. Then
$I_{\theta(x_s)}=\{s\}\cup I_v$.  Since $x_s$ is normal,
so is $\theta(x_s)$.  By Lemma \ref{xxlem2.5}(a), $I_{\theta(x_s)}
\subset B_w$, and whence $s\in B_w$, a contradiction.
Therefore condition (i) holds.
Let $g$ be an elementary transformation such that $g(x_i)=v$
for some $i\in B_w$ and $g(x_s)$ is fixed by $\theta$ for all
$s\in B_w\setminus \{i\}$. Then $g^{-1}\theta g=\theta_{i,\lambda}$.
\end{proof}

A partition ${\mathcal D}=\{D_w\;|\; w\in W\}$
of $[n]$ (that is, $[n]$ is a disjoint union
$\bigcup_{w\in W}D_w$) is called a $\pp$-partition
if, for any two distinct $i,i'\in D_w$,
$p_{ij}=p_{i'j}$ for all $j\in [n]\setminus D_w$.
The block decomposition is a $\pp$-partition.
For a given partition ${\mathcal D}$, let $\Aut_w(A)$
be the subgroup of $\Aut(A)$ consisting of the
automorphisms $g$ satisfying
$$\begin{cases} g(x_i)=x_i& i\not\in D_w, \quad {\text{and}} \\
g(x_i) \in \sum_{s\in D_w} kx_s & i\in D_w.\end{cases}$$

\begin{lemma}
\label{xxlem2.7} Fix a $\pp$-partition ${\mathcal D}
=\{D_w\;|\; w\in W\}$.
Let $G$ be a finite group of graded algebra
automorphisms of $A$. Let $G_w=G\cap \Aut_w(A)$.
In parts (d,e) suppose that $G=\prod_{w\in W} G_w$ .
\begin{enumerate}
\item
For any $w$, pick any $j\in D_w$ and define
$$\phi_{w}: x_i\to \begin{cases}
q_{ji}x_i & i\not\in D_w\\
x_i& i\in D_w.\end{cases}$$
Then $\phi_w$ commutes with automorphisms in
$\Aut_{w'}(A)$ for all $w'$.
\item
Let $A_w$ be the subalgebra generated by
$\sum_{i\in D_w} kx_i$. Let $\phi=\{\phi_w\}$ as in part (a).
We can use this set of automorphisms to twist $A$; denote the
resulting algebra by  $A^\phi$.
Then $A^\phi$ is isomorphic to $\bigotimes_{w\in W} A_w$, and
hence $A$ is a graded twist of the tensor product
$\bigotimes_{w\in W} A_w$.
\item
If each $D_w$ is a subset of some $B_{w'}$, then
$A^\phi$ is a commutative polynomial ring.
\item
The group $G$ acts naturally on the twisted
algebra $A^\phi$, and $(A^\phi)^G$ is isomorphic to a
twisted algebra $(A^G)^\phi$, where $\phi$ is a set of
graded automorphisms of $A^G$ induced from $A$.
\item
The fixed ring $A^G$ has finite global dimension
(respectively, is regular) if and only if $(A^G)^\phi$ has
finite global dimension (respectively, is regular) if and only if
each part $A_w^{G_w}$ has finite global dimension
(respectively, is regular).
\end{enumerate}
\end{lemma}

\begin{proof}
(a) Straightforward. Recall that $q_{ij}=\sqrt{p_{ij}}$.

(b) Use the proof of Lemma \ref{xxlem2.2}(e)

(c) This follows from part (b) and the fact that
each $A_w$ is a commutative polynomial ring.

(d) Since any $\phi_w$ commutes with elements
in $G$, $G$ acts on $A^\phi$ naturally. The rest is
easy to check.

(e) The first assertion is Lemma \ref{xxlem1.7}(a).
By part (d) we may assume $\phi$ is trivial, and hence
$A$ is a tensor product of $A_w$'s and $A^G$ is a
tensor product of $A_w^{G_w}$. Then assertion follows
from Lemma \ref{xxlem1.6}.
\end{proof}

\begin{lemma}
\label{xxlem2.8} Fix a $\pp$-partition ${\mathcal D}
=\{D_w\;|\; w\in W\}$ and let $\phi_w$
be as in Lemma \ref{xxlem2.7}(a). Suppose $G=\prod_{w\in W} G_w$.
Then an element $\eta$ in $G$ is a reflection
(respectively, mystic reflection) of $A$ if and only if
it is a reflection (respectively, mystic reflection) of $A^\phi$.
\end{lemma}

\begin{proof} By Lemma \ref{xxlem2.7}(a) and by
the hypothesis that $G=\prod_{w\in W} G_w$, $\phi_w$
commutes with every element in $G$. Therefore every
element in $G$ is also a graded algebra automorphism
of $A^\phi$ when $A^\phi$ is identified with $A$ as
graded vector spaces. This also implies that
$Tr_A(\eta,t)=Tr_{A^\phi}(\eta,t)$ for all $\eta\in G$.
Therefore $\eta$ is a quasi-reflection of $A$ if
and only if it is a quasi-reflection of $A^\phi$.
Since $\det \eta|_{A_1}=\det \eta|_{A^\phi_1}$,
$\eta$ is a mystic reflection of $A$ if and only if
it is a mystic reflection of $A^{\phi}$ by Proposition
\ref{xxprop1.4}(b).
\end{proof}

%
%


\section{Circles induced by finite group actions}
\label{sec3}

In this section we will study further the structure of
$A$ when $A$ admits a mystic reflection. As in the last section
we fix a parameter set $\pp$. Our starting
point is the following proposition.

\begin{proposition}
\label{xxprop3.1}
Let $\tau$ be a mystic reflection of $A$.
\begin{enumerate}
\item
Up to a sequence of elementary transformations,
$\tau$ is a standard mystic reflection $\tau_{i,j,1}$
for some $i,j\in [n]$.
\item
There is a pair $(i,j)$ such that $p_{ij}=-1$.
\item
Let $i,j\in [n]$ be the integers described in part (a).
Then $B(i)=\{i\}$ and $B(j)=\{j\}$. As a consequence,
$\tau=\tau_{i,j,\lambda}$ (without conjugating elementary
transformations).
\end{enumerate}
\end{proposition}

\begin{proof}
By Proposition \ref{xxprop1.4}(b,c), there is a basis
$\{y_i\;|\; i\in [n]\}$ such that
\begin{equation}
\label{xxequ3.1.1}
\tau(y_i)=\begin{cases}
y_i & i\leq n-2\\
y_n & i=n-1\\
-y_{n-1}& i=n
\end{cases}
\tag{3.1.1}
\end{equation}
and $y_ny_{n-1}=-y_{n-1}y_n$. Since $y_ny_{n-1}=-y_{n-1}y_n$,
Lemma \ref{xxlem2.5}(b) implies that $I_{y_n}\cap I_{y_{n-1}}
=\emptyset$,  $p_{ij}=-1$ for all $i\in I_{y_n}$, and
$j \in I_{y_{n-1}}$. Thus we have proved part (b).

We claim that part (c) follows from part (a). By part (a)
and induction on the number of elementary transformations
we may assume that $\tau= g \tau_{i,j,1} g^{-1}$ where $g$ is
an elementary transformation. Suppose $|B(i)|=\alpha>1$ and then
write $B(i)=\{i_1=i,i_2,\cdots i_\alpha\}$. We have $g(x_{i_s})=
x'_{i_s}\in \bigoplus_{t=1}^\alpha kx_{i_t}$ for all
$s=1,\cdots,\alpha$. Similarly, we write $B(j)=\{j_1=j,
\cdots, j_{\beta}\}$ and $g(x_{j_s})=x'_{j_s}\in
\bigoplus_{t=1}^\beta kx_{j_t}$ for all
$s=1,\cdots,\beta$. Since $\tau= g \tau_{i,j,1} g^{-1}$,
we have $\tau(x'_i)=x'_j$, $\tau(x'_j)=-x'_i$,
and $\tau(x'_s)=x'_s$ for all $s\neq i,j$. By the last
paragraph, $p_{i_sj_t}=-1$ for some $i_s\in B(i)$ and
$j_t\in B(j)$. By the definition of the blocks,
we have $p_{i_sj_t}=-1$ for all $i_s\in B(i)$ and
$j_t\in B(j)$. Since $x'_i=x'_{i_1}$
commutes with $x'_{i_\alpha}$ (they are in the same block),
$x'_j=x'_{j_1}$ commutes with $x'_{i_\alpha}$ after applying
$\tau$. This contradicts the fact that $p_{i_sj_t}=-1$.
Therefore $B(i)=\{i\}$. By symmetry, $B(j)=\{j\}$.
Up to a scalar, we may assume that $x'_i=x_i$
and $x'_j=x_j$. Hence $\tau(x_i)=x_j$, $\tau(x_j)=-x_i$,
and $\tau(x'_s)=x'_s$ for all $s\neq i,j$. For $x_s$ outside
of the blocks $B(i)$ and $B(j)$, we may assume
$g$ is the identity. Therefore $\tau=\tau_{i,j,1}$ up to
a scalar change on $x_i$.

It remains to prove (a). We assume \eqref{xxequ3.1.1}
and, without loss of generality, we may assume that, up to
elementary transformations, $|I_{y_{n-1}}|$ is minimal.

Case 1: $|I_{y_{n-1}}|=1$, or, equivalently, up to a scalar,
$y_{n-1}=x_t$ for some $t\in [n]$. Up to a
permutation we may assume that $y_{n-1}=x_{n-1}$. Then
$y_{n-1}$ is normal and so is $y_n=\tau(y_{n-1})$.
Lemma \ref{xxlem2.5}(a) says that $I_{y_n}$ is a subset
of a block $B_w$ for some $w$ (denote this block as
$B(I_{y_n})$). Using an elementary transformation and a
permutation if necessary, we may assume that $y_n=x_n$.
Since $\tau$ becomes the identity on $A/(x_{n-1},x_n)$,
we have
$$\tau(x_i)=x_i+a_i x_{n-1}+b_i x_n$$
for all $i\in [n-2]$. Recall that $p_{n-1n}=-1$.
Using Lemma \ref{xxlem2.5}(a) we see the following:
\begin{enumerate}
\item[(i)]
If $i\in B(n)$, then $a_i=0$ because $n-1\not\in B(n)$. Now
the fact that $\tau^2(x_i)=x_i+b_i x_n-b_ix_{n-1}$ is normal
implies that $b_i=0$.
\item[(ii)]
If $i\in B(n-1)$, then $b_i=0$. A similar argument to that above
shows that $a_i=0$.
\item[(iii)]
If $i\not\in B(n)\cup B(n-1)$, then $a_i=b_i=0$.
\end{enumerate}
Therefore $\tau=\tau_{n-1,n,1}$ up to an elementary
transformation and a permutation.

Case 2: $|I_{y_n}|=1$. This is similar to Case 1.

Case 3: $|I_{y_{n-1}}|>1$ and $|I_{y_n}|>1$.
Since $\tau$ is the identity on $A_1/(ky_{n-1}+ky_{n})$,
we have
$$\tau(x_i)=x_i+a_i y_{n-1}+b_iy_n$$
for all $i\in [n]$. Since $x_i$ is normal, so is $\tau(x_i)$.
Let $i\not\in I_{y_n}\cup I_{y_{n-1}}$. If $a_i\neq 0$,
then $B(I_{\tau(x_i)})\supset \{i\}\cup I_{y_{n-1}}$ by
Lemma \ref{xxlem2.5}(a). This implies that $I_{y_{n-1}}$ is a single
element after an elementary transformation, a contradiction.
Therefore $a_i=0$. Similarly, $b_i=0$ for all  $i\not\in
I_{y_n}\cup I_{y_{n-1}}$.

For the final argument we consider two cases.

First we consider the case when $|I_{y_{n}}|=|I_{y_{n-1}}|=2$.
Without loss of generality, we may assume that
$y_{n-1}=x_1+x_2$ and $y_n=x_{n-1}+x_n$. So we have
$\tau(x_i)=x_i$ for all $2<i<n-1$. So we have
$$\begin{aligned}
\tau(x_1)&=x_1+a_1(x_1+x_2)+b_1(x_{n-1}+x_n), \text{ and}\\
\tau(x_2)&=x_2+a_2(x_1+x_2)+b_2(x_{n-1}+x_n).\\
\end{aligned}
$$
Since $\tau(x_1+x_2)=\tau(y_{n-1})=y_n=x_{n-1}+x_n$,
either $b_1$ or $b_2$ is nonzero. Suppose $b_1\neq 0$.
Then $\{n-1,n\}\in B(I_{\tau(x_1)})$. Therefore up to
an elementary transformation, $|I_{y_n}|=1$.
But this was done in Case 1.

Second we assume that $|I_{y_{n-1}}|\geq 3$.
Write
$$\tau(x_i)=x_i+a_i y_{n-1}+b_i y_n.$$
Since $I_{y_{n-1}}\cap I_{y_n}=\emptyset$,
$I_{\tau(x_i)}\supset I_{y_{n-1}}\setminus \{i\}$
if $a_i\neq 0$. This means that we can strictly reduce
the cardinality of $I_{y_{n-1}}$ by an elementary
transformation. This contradicts the minimality of
$|I_{y_{n-1}}|$. Therefore $a_i=0$ for all $i$. Thus
$\tau$ is an identity on $A_1/(ky_{n})$, which yields
a contradiction. By symmetry, it is impossible to have
$|I_{y_{n}}|\geq 3$. This completes the proof.
\end{proof}

%


\begin{definition}
\label{xxdefn3.2}
Fix a finite subgroup $G$ of $\Aut(A)$.
\begin{enumerate}
\item
Any single element subset of $[n]$ is called a {\it trivial circle}.
\item
Let $C$ be a subset of $[n]$ consisting of at least two integers.
We call $C$ a {\it circle} if, for each pair of distinct integers
$(i, j)\subset C$, there is a sequence of mystic reflections
$\{\tau_{i_s,j_s,\lambda_s}\in G\;|\; s=0,\cdots , t\}$ for some
$\lambda_s\in k^{\times}$ such that $i=i_0, i_{s+1}=j_s$ for all $s$,
and $j_{t}=j$. By Proposition \ref{xxprop3.1}(c), every non-trivial
circle is a union of some single-element blocks.
\item
We call $C$ a {\it maximal circle} if $C$ is a non-trivial
circle and it is not properly contained in another circle.
\item
A {\it circle decomposition} is a disjoint union of
maximal circles and trivial circles
$$[n]=\bigcup_{u\in U} C_u.$$
\item
A {\it block-circle decomposition} is a disjoint union
$$[n]=\bigcup_{u\in U_1} C_u \cup \bigcup_{w\in W_1}
B_w$$
where each $C_u$ for $u\in U_1$ is a
maximal circle and where each $B_w$ for $w\in W_1$
is a block in $[n]$. This decomposition can be formed
by first having a block decomposition and then
joining single-element blocks together to possible
maximal circles.

For simplicity, we write block-circle decomposition as
$$[n]=\bigcup_{v\in V} D_v$$
where $V=U_1\cup W_1$ and $D_v$ is either $B_v$ or $C_v$.
\end{enumerate}
\end{definition}

The block-circle decomposition always exists and is
uniquely determined by the parameters $(\pp,G)$.
As we noted before, the block decomposition is a
$\pp$-partition. Next we show that the block-circle
decomposition is a $\pp$-partition.

\begin{proposition}
\label{xxprop3.3}
Let $G$ be a finite group of $\Aut(A)$ and let
$[n]=\bigcup_{v\in V} D_v$ be the block-circle decomposition.
Let $A_v$ be the subalgebra generated by $\sum_{i\in D_v}k x_i$.
\begin{enumerate}
\item
Let $G_v=G\cap \Aut_v(A)$. Assume that $G$ is generated by
quasi-reflections. Then $G=\prod_{v\in V} G_v$ and each
$G_v$ is generated by quasi-reflections.
\item
The block-circle decomposition is a $\pp$-partition, namely,
for any $i\neq i'$ in the same $D_v$ and $j\not\in D_v$,
$p_{ij}=p_{i'j}$.
\end{enumerate}
\end{proposition}

\begin{proof} (a) Note that the block decomposition is a
sub-partition of the circle-block decomposition.
By Proposition \ref{xxprop2.6}, every reflection  is an elementary
transformation, hence it belongs to $\Aut_w(A) \cap G = G_w$ for some $w$. By Proposition
\ref{xxprop3.1} every mystic reflection is in $G_v$ for some
$v\in U_1$. Therefore the assertion follows.

(b) We need to show this only for $v\in U_1$. By induction we
may assume that there is a mystic reflection of the form
$\tau_{i,i',1}$. Since $\tau_{i,i',1}(x_j)=x_j$ for
all $j\not\in D_v$, $\tau_{i,i',1}(x_i)=x_{i'}$ implies
that $p_{ij}=p_{i'j}$.
\end{proof}

The next corollary follows from Proposition \ref{xxprop3.3},
Lemma \ref{xxlem2.7} and Lemma \ref{xxlem4.1} in the next
section.

\begin{corollary}
\label{xxcor3.4}
Let $G$ be a finite group of $\Aut(A)$ and let
$[n]=\bigcup_{v\in V} D_v$ be the block-circle decomposition.
Let $A_v$ be the subalgebra generated by $\sum_{i\in D_v}k x_i$.
\begin{enumerate}
\item
For any $v$ and any $j\in D_v$, there is an algebra
automorphism determined by
$$\phi_{v}: x_i\to \begin{cases}
q_{ji}x_i & i\not\in D_v\\
x_i& i\in D_v.\end{cases}$$
Further $\phi_v$ commutes with automorphisms in
$\Aut_{v'}(A)$ for all $v'$.
\item
Let $\phi=\{\phi_v\;|\; v\in V\}$. We can
use this set of automorphisms to twist $A$; denote the
resulting algebra by $A^\phi$.
Then $A^\phi$ is isomorphic to $\bigotimes_{v\in V} A_v$, and
hence $A$ is a graded twist of
$\bigotimes_{v\in V} A_v$.
\item
The fixed ring $A^G$ has finite global dimension
(respectively, is regular) if and only if $(A^G)^\phi$ has
finite global dimension (respectively, is regular) if and only if
each part $A_v^{G_v}$ has finite global dimension
(respectively, is regular).
\item
If $D_v$ is a block, then $A_v$ is a commutative polynomial
ring.
\item
If $D_v$ is a circle, then $A_v$ is isomorphic to
$k_{-1}[x_1,\cdots,x_m]$ where $m=|D_v|$.
\item
If $G$ is generated by quasi-reflections of $A$, then
each $G_v$, when restricted to $A_v$, is generated by
quasi-reflections of $A_v$.
\end{enumerate}
\end{corollary}

\begin{proof} By Proposition \ref{xxprop3.3},
the block-circle decomposition is a $\pp$-partition.
See Lemma \ref{xxlem2.7} for parts (a,b,c,d).

Part (e) will be proved in Lemma \ref{xxlem4.1} in the next
section.

(f) By parts (b,c) we may assume that $A$ is the tensor
product $\bigotimes_{v\in V} A_v$. Then the trace formula
implies that $g\in G_v$ is a quasi-reflection if and only
if $g|_{A_v}$ is a quasi-reflection. The assertion
follows because $G_v$ is generated by quasi-reflections
of $A$ by the proof of Proposition \ref{xxprop3.3}(a).
\end{proof}

By Corollary \ref{xxcor3.4}(c,f) we need to
prove the Shephard-Todd-Chevalley Theorem only for
each block/circle  $A_v$. If $D_v$ is a block,
then the classical Shephard-Todd-Chevalley Theorem
applies. The next section deals with the circle case.

\section{One circle case and the proof of the Main Theorem}
\label{sec4}

In most of this section we
let $G$ be a finite subgroup of $\Aut(A)$ that contains
mystic reflections such that $[n]$ is a circle. Assume
that $n\geq 2$.

\begin{lemma}
\label{xxlem4.1} Suppose $[n]$ is a circle for the group
$G$, which is generated by quasi-reflections of $A$.
Up to scalar change of a basis for $A_1$, $\tau_{i,j,1} \in G$
for all $i\neq j$. As a consequence, $p_{ij}=-1$ for all $i<j$.
\end{lemma}

\begin{proof} We use induction on $m$ to show that
$\tau_{i,j,1}\in G$ for all $1\leq i,j\leq m$. Nothing
is to be proved for $m=1$, so we begin with $m=2$. By
Definition \ref{xxdefn3.2}(b), there is a sequence of
$\{\tau_{i_s,j_s,\lambda_s}\in G\;|\;
s=0, \cdots, t\}$ such that $1=i_0,
i_{s+1}=j_s, j_t=2$. Choose $t$ to be minimal; then
we claim that $t=0$. If $t>0$, we have
$\tau^{-1}_{i_1,j_1,\lambda_1} = \tau_{i_1,j_1,-\lambda_1}$, and
$$\tau_{1,j_1,\lambda_0\lambda_1}=
\tau_{i_1,j_1,\lambda_1}
\tau_{1,i_1,\lambda_0}
\tau_{i_1,j_1,\lambda_1}^{-1}.$$
Then we can make a shorter sequence, a contradiction.
Therefore $t=0$ and $\tau_{1,2,\lambda}\in G$ for some
$\lambda\in k^{\times}$. Replacing $x_2$ with $\lambda x_2$,
we have $\tau_{1,2,1}\in G$, and this completes the
proof for $m=2$. Suppose now $m>2$ and the assertion
holds for $m-1$. By a similar argument, we have
$\tau_{m-1,m,\lambda}\in G$. Letting the new $x_m$ be
$\lambda x_m$, we have $\tau_{m-1,m,1}\in G$. For any
$i<m-1$, we have
$$\tau_{i,m,1}=
\tau_{m-1,m,1}
\tau_{i,m-1,1}
\tau_{m-1,m,1}^{-1}\in G.$$
By the induction hypothesis, $\tau_{i,j,1}\in G$ for all
$1\leq i,j \leq m$. The assertion follows now by
induction on $m$.
\end{proof}

We now assume that $p_{ij}=-1$ for all $i<j$; fix a
basis $\{x_i\;|\; i\in [n]\}$ so that $\tau_{i,j,1}\in G$
for all $i\neq j$. Let
$$\Theta_i=\{\lambda\in k^{\times}\; |\; \theta_{i,\lambda}\in G\},$$
$$T_{i,j}=\{\lambda\in k^{\times}\;|\;
\tau_{i,j,\lambda}\in G\},$$
and
$$S_{i,j}=\{\lambda\in k^{\times}\;|\;
s_{i,j,\lambda}\in G\},$$
where $s_{i,j,\lambda}$ is the automorphism determined by
$$s_{i,j,\lambda}(x_s)=\begin{cases}
x_s& s\neq i,j\\
\lambda x_i &s=i\\
\lambda^{-1}x_j& s=j.\end{cases}$$

\begin{lemma}
\label{xxlem4.2} Suppose $[n]$ is a circle for a group
$G$ that is generated by quasi-reflections of $A$.
\begin{enumerate}
\item
$\Theta_i=\Theta_j$ for all $i,j$, and $\Theta_i\cong
{\mathbb Z}/(\alpha)$ for some integer $\alpha$.
\item
$S_{i,j}=S_{i',j'}$ for all $(i,j)$, and $(i',j')$
and $S_{i,j}\cong
{\mathbb Z}/(\beta)$ for some even integer $\beta$.
\item
$T_{i,j}=T_{i',j'}$ for all $(i,j)$ and $(i',j')$,
and $T_{i,j}=S_{i,j}\cong
{\mathbb Z}/(\beta)$ for some even integer $\beta$.
\item
$\alpha$ divides $\beta$.
\end{enumerate}
\end{lemma}

\begin{proof}
First note that $\Theta_i, T_{i,j}$, and $S_{i,j}$
are all subgroups of $k^{\times}$.

(a) The first assertion follows from
$$\theta_{j,\lambda}=\tau_{i,j,1}\theta_{i,\lambda}
\tau_{i,j,1}^{-1}.$$
The second follows from the fact that any finite subgroup
of $k^{\times}$ is cyclic, whence it is
of the form $\{\lambda\;|\; \lambda^\alpha=1\}$
for some $\alpha$.

(b) Since $S_{i,j}=S_{j,i}$ for $i\neq j$, we may assume
that $i=i'$ without loss of generality. In this case $i,j,j'$
are all different. The first assertion follows from the fact that
$$s_{i,j',\lambda}=\tau_{j',j,1}s_{i,j,\lambda}
\tau_{j',j,1}^{-1}.$$
Similar to the proof of (a), $S_{i,j}\cong {\mathbb Z}/(\beta)$.
The reason that $\beta$ is even is that $s_{i,j,-1}=\tau_{i,j,1}^2$
is in $G$.

(c) This is true because $\tau_{i,j,1}^{-1}\tau_{i,j,\lambda}=
s_{i,j,\lambda}$.

(d) This is true because $s_{i,j,\lambda}=\theta_{i,\lambda}
\theta^{-1}_{j,\lambda}$.
\end{proof}

Suppose $\alpha|\beta$ and $2|\beta$.
Let $M({n,\alpha,\beta})$ be the subgroup of $\Aut(A)$ generated
by $\{\theta_{i,\lambda}\;|\; \lambda^\alpha=1\}\cup
\{\tau_{i,j,\lambda}\;|\; \lambda^{\beta}=1\}$.
Lemma \ref{xxlem4.2}, Proposition \ref{xxprop2.6}, and Proposition
\ref{xxprop3.1}(c) imply the following fact. Recall
that in this section (except for the last theorem)
$A=k_{-1}[x_1,\cdots,x_n]$.

\begin{lemma}
\label{xxlem4.3}
If $G$ is a finite subgroup of $\Aut(A)$ generated by
quasi-reflections and if $[n]$ is circle, then $G$ is
$M({n,\alpha,\beta})$ for some $\alpha$ and $\beta$.
\end{lemma}

Let $G=M({n,\alpha,\beta})$.
Let $B$ be the subalgebra of $A$ generated by $x_1^{\alpha},
\cdots, x_n^{\alpha}$. If $\alpha$ is even then $B$ is
the commutative polynomial ring. If $\alpha$ is odd, then
$B\cong A$ under the map that sends $x_i\to x_i^\alpha$. Let $z_i=x_i^\alpha$
for all $i$.
To cover both cases of $\alpha$, we write $B=k_{\pm 1}[z_1,
\cdots,z_n]$. Let $R$ be the subgroup generated by
all $\theta_{i,\lambda}$. Then $B=A^{R}$. Let $G_B$
be the induced group of $G$ acting on $B$.  Note the
similarity between the invariants described below and those
obtained when the symmetric group $S_n$ acts on
$k[x_1, \cdots, x_n]$.

\begin{proposition}
\label{xxprop4.4} Retain the notation described above.
\begin{enumerate}
\item
$A^G=B^{G_B}$ as subalgebras of $A$.
\item
If $\alpha$ is even, then $G_B$ is a classical
reflection group of the commutative polynomial ring
$B$. As  a consequence, $B^{G_B}$ is regular.
\item
If $\alpha$ is odd, then $G_B$ contains no reflections
of $B$. In this case $G_B$ is generated by mystic reflections
of $B$.
\item
Suppose $\alpha$ is odd. Let $2t$ be the $\beta$ for the
group $G_B$. Then $B^{G_B}$ is the commutative polynomial ring
$D:=k[(z_1\cdots z_n),\sum_i z_i^{2t}, z_i^{4t},\cdots
\sum_i z_i^{(n-1)2t}]$. As a consequence, $B^{G_B}$ is regular.
\end{enumerate}
By (b,d), $A^G$ is always isomorphic to a
commutative polynomial ring.
\end{proposition}

\begin{proof} (a) Let $R$ be the subgroup generated by the
reflections $\theta_{i,\lambda}$. Then $R$ is a normal
subgroup of $G$, and it is easy to see that $A^R=B$. Since
$R$ is normal, $G$ induces a subgroup $G_B$ of $\Aut(B)$
that is generated by the induced maps of $\tau_{i,j,\lambda}$,
denoted by $\tau'_{i,j,\lambda}$. Clearly, $A^G=B^{G_B}$.

(b) If $\alpha$ is even, all of the $\tau'_{i,j,\lambda}$
are reflections of $B$. Hence $B^{G_B}$ is regular.

(c) If $\alpha$ is odd, then all of the $\tau'_{i,j,\lambda}$
are equal to $\tau_{i,j,\lambda^\alpha}$ for the skew polynomial
ring $B=k_{-1}[z_1,\cdots,z_n]$, and so these are mystic reflections.
By Proposition~\ref{xxprop1.4}(d), $G_B$ does not contain any reflection.

(d) Let $C=k[((z_1\cdots z_n)^{2t},\sum_i z_i^{2t},\cdots,
\sum_i z_i^{(2n-2)t}]$. It is clear that
$C$ is a subalgebra of $B^{G_B}\cap k[z_1^{2t},
\cdots,z_n^{2t}]$. For any homogeneous element $w\in B$,
we define
$$\Gamma(w):={\frac{1}{|G_B|}}
\sum_{g\in G_B} g(w).$$
It suffices to show that $\Gamma(w)\in D$ for all $w\in B$.
(The definition of $D$ is given in Proposition \ref{xxprop4.4}(d).)
Since $\Gamma$ is linear, we may take $w$ to be a
monomial. The following analysis uses induction on the degree
of $w$.

Case 1: Suppose $(z_1\cdots z_n)| w$. Then $w=
(z_1\cdots z_n) w'$ for some $w'$ with its degree less than
the degree of $w$, and $\Gamma(w)=(z_1\cdots z_n)\Gamma(w')$.
The assertion follows by induction.

Case 2: Suppose $(z_1\cdots z_n)\nmid w$.  If $\Gamma(w)=0$
we are done.
Otherwise $\Gamma(w)=\sum_{(d)} c_{(d)}
z_1^{d_1}\cdots z_n^{d_n}\neq 0$ where
$(d)=(d_1,\cdots,d_n)$. By the definition of
$\Gamma(w)$ and the hypothesis that $(z_1\cdots z_n)\nmid w$,
we see that for each $c_{(d)}\neq 0$
there is an $i$ such that $d_i=0$. If, for each $c_{(d)}\neq 0$,
$2t\mid d_i$ for all $i$, then $\Gamma(w)$ is in
\[k[z_1^{2t}, \cdots, z_n^{2t}]^{G_B},\]
which is the ring of symmetric polynomials in
$\{z_1^{2t}, \ldots, s_n^{2t}\}$.
By the Newton identities (see e.g. \cite[p. 317]{CLO}) this is
$k[(z_1^{2t} \cdots z_n^{2t}),\sum_i z_i^{2t},\cdots,
\sum_i z_i^{(2n-2)t}] = C \subseteq D$, and
we are done in this case. Otherwise there is
$c_{(d)}\neq 0$ such that $2t\nmid d_i$ for some $i$.
Without loss of
generality, we may assume that:
\begin{enumerate}
\item
$c_{(d)}\neq 0$ for some $(d)=(d_1,\cdots,d_{n-1},0)$,
(this means that we pick $d_n=0$),
\item
$2t\nmid d_i$ for some $i$, and
\item
$2t| d_j$ for all $j>i$.
\end{enumerate}
Note that every $\tau_{i,j,\lambda}$ maps a monomial in
$\Gamma(w)$ to another monomial in $\Gamma(w)$ (possibly
the same). In particular,
$\tau_{i,n,\lambda}$ maps $c_{(d)}z_1^{d_1}\cdots z_i^{d_i}
\cdots z_{n-1}^{d_{n-1}}$ to
$$c_{(d)}z_1^{d_1}\cdots (\lambda z_n)^{d_i}
\cdots z_{n-1}^{d_{n-1}}=\lambda^{d_i}c_{(d)}z_1^{d_1}\cdots \widehat{z_i}
\cdots z_{n-1}^{d_{n-1}} z_n^{d_i}.$$
Since $\tau_{i,n,\lambda}(\Gamma(w))=\Gamma(w)=\tau_{i,n,1}
(\Gamma(w))$,
the element $\lambda^{d_i}c_{(d)}z_1^{d_1}\cdots \widehat{z_i}
\cdots z_{n-1}^{d_{n-1}} z_n^{d_i}$ is independent of
$\lambda$. This implies that $\lambda^{d_i}=1$ or
$2t\mid d_i$, a contradiction.
This completes the proof.
\end{proof}

Now we are ready to prove our main result, Theorem \ref{xxthm0.1}.

\begin{theorem}
\label{xxthm4.5} Let $A$ be a skew polynomial ring
$k_{p_{ij}}[x_1,\cdots, x_n]$ and
$G$ a finite subgroup of $\Aut(A)$. Then the fixed subring
$A^G$ has finite global dimension if and only if $G$ is
generated by quasi-reflections. In this case the fixed subring
$A^G$ is isomorphic to a skew polynomial ring with
possibly different $\{p_{ij}\}$.
\end{theorem}

\begin{proof} By Proposition 1.5 we need to show only one
direction, namely, that if $G$ is generated by
quasi-reflections then $A^G$ has finite global dimension.
It suffices to show that $A^G$ is isomorphic to
a skew polynomial ring.

By Corollary \ref{xxcor3.4}(c,f) and Lemma \ref{xxlem1.7}(b),
we need to show the claim for $A_v$, where $D_v$ is either a
block or a circle in the block-circle
decomposition of $[n]$. If $D_v$ is a block, then $A_v$ is a
commutative polynomial ring, and the claim follows from
the classical Shephard-Todd-Chevalley Theorem. If
$D_v$ is a circle, the claim follows from Lemma
\ref{xxlem4.3} and Proposition \ref{xxprop4.4}.
\end{proof}

\begin{corollary}
\label{xxcor4.6}  Let $A=k_{p_{ij}}[x_1,\cdots, x_n]$ and
let $G$ be a finite subgroup of $\Aut(A)$. If $G$ is
generated by quasi-reflections, then, as an abstract
group, $G$ is a product of classical reflection groups and
copies of some  $M(n,\alpha,\beta)$s.
\end{corollary}

\begin{proof} Using the block-circle decomposition
$[n]=\bigcup_{v\in V} D_v$ [Definition \ref{xxdefn3.2}(e)], 
$G$ is a product $\prod_{v\in V} G_v$ and each $G_v$ is generated by
quasi-reflections [Proposition \ref{xxprop3.3}]. If
$D_v$ is a block, then $G_v$ is isomorphic to a classical
reflection group [Corollary \ref{xxcor3.4}(d)]. If $D_v$
is a circle, then $G_v$ is isomorphic to $M(n,\alpha,\beta)$
by Lemma \ref{xxlem4.3}. The assertion follows.
\end{proof}

\begin{remark}
\label{xxrem4.7}
We will see that there are ``mystic reflection groups''
$M(n,\alpha,\beta)$ that cannot be realized
as classical reflection groups 
[Example \ref{xxex6.2} and Example \ref{xxex6.3}].
\end{remark}

The following conjecture, which is suggested by Corollary
\ref{xxcor4.6}, is related to Conjecture \ref{xxcon0.2}.

\begin{conjecture}
\label{xxcon4.8}
Let $B$ be a quantum polynomial ring. Suppose $G$ is a finite
subgroup of $\Aut(B)$ such that $B^G$ has finite global dimension.
Then $G$ is a product of classical reflection groups and
copies of $M(n,\alpha,\beta)$s.
\end{conjecture}

\section{Toward a quantum Shephard-Todd-Chevalley theorem}
\label{sec5}

In the first half of this section we prove Theorem \ref{xxthm0.3},
namely, a version of Shephard-Todd-Chevalley Theorem for general
quantum polynomial rings when the group is abelian. The proof of
Theorem \ref{xxthm0.3} is quite different from the proof of
Theorem \ref{xxthm0.1}.

\begin{lemma}
\label{xxlem5.1}
Let $B$ be a quantum polynomial ring and let $G$ be a finite
abelian subgroup of $\Aut(B.)$  Suppose that
$\{g_1, g_2, \ldots , g_r\}$ is a minimal generating set
for $G$ where each $g_i$ is a quasi-reflection.
If $i \neq j$ and $y \in B_1$ is a common eigenvector
for $g_i$ and $g_j$ with $g_i(y) = \lambda_i y$ and $g_j(y)
= \lambda_j y$, then either $\lambda_i =1$ or $\lambda_j =1$.
\end{lemma}

\begin{proof}
Since $G$ is abelian, we can find a $k$-linear basis
$\{y_1, y_2, \ldots, y_n\}$ of $B_1$ such that the action
of each $g_i$ on $B_1$ is diagonal with respect to this
basis.  Suppose that neither $\lambda_i$ nor $\lambda_j$ is
equal to $1$. For simplicity  assume that $g_i=g_1$ and
that $g_j=g_2.$ Suppose to the contrary that $\lambda_1\neq 1$
and $\lambda_2\neq 1$. We will obtain a contradiction
in each of the following four cases.

Case 1: Suppose that $g_1$ and $g_2$ are both reflections.
Since the dimensions of the corresponding non-invariant
eigenspaces are each $1$ [Definition \ref{xxdefn1.3}(a)],
there is no loss of generality in assuming that $y= y_n$.
By Proposition \ref{xxprop1.4}(a)  $y_n$ is
a normal element of $B$. In this case $\lambda_1$ and
$\lambda_2$ are $m_1$th and $m_2$th roots of unity,
respectively. Let $\lambda$ be the generator of the subgroup
$\langle \lambda_1, \lambda_2 \rangle$ of $k^{\times}$. Then
$\lambda$ is an $m$th root of unity for $m =\rm{lcm}(m_1, m_2)$,
and $\lambda = \lambda_1^{r_1} \lambda_2^{r_2}$ for integers
$r_1$ and $r_2$.  Let $g = g_1^{r_1}g_2^{r_2}$.  Then $g$ is
an automorphism of $B$ with $g(y_n) = \lambda y_n$ and
$g(y_s)=y_s$ for $s \neq n$.  Thus $g$ induces the identity
automorphism on the factor algebra $\bar{B} =B/(y_n)$ so
that
$$Tr_{\bar{B}}(g,t) = \frac{1}{(1-t)^{n-1}}.$$
But it is also the case that $Tr_{\bar{B}}(g,t) =
(1-\lambda t)Tr_B(g,t).$  Hence
$$ Tr_B(g,t) = \frac{1}{(1-t)^{n-1}(1-\lambda t)},$$
and $g$ is a quasi-reflection.  Since both $g_1$ and $g_2$ are
powers of $g$, we have contradicted the minimality of $\{g_1, g_2, \ldots , g_r\}$ as a generating
set.

Case 2:  Suppose that $g_1$ is a reflection and $g_2$ is a
mystic reflection.  Then as in Case 1 we may assume that $y=y_n$,
and, since $g_1$ is a reflection, that $y_n$ is a normal element of
$B$.  By Definition \ref{xxdefn1.3}(b) we  have that $g_2(y_n) =
\pm i y_n,$; without loss of generality, say $g_2(y_n) = iy_n$.
There then must be another eigenvector, which we may take as
$y_{n-1}$, such that $g_2(y_{n-1})=-iy_{n-1}$. Furthermore, we
assume that the subalgebra generated by $y_{n-1}$ and $y_n$ is
isomorphic to $C:=k\langle y_{n-1}, y_n \rangle/(y_{n-1}^2-y_n^2)$
[Proposition \ref{xxprop1.4}(c)].  Since $y_n$ is a normal
element of $B$ we have that
\[ y_n y_{n-1} = (a_{n-1}y_{n-1} +a_ny_n + \sum_{i=1}^{n-2}a_iy_i)y_n. \]
Applying $g_2$ we have
\[ y_n y_{n-1} = (a_{n-1}y_{n-1} -a_ny_n + i \sum_{i=1}^{n-2}a_iy_i)y_n. \]
Subtracting gives
\[ (2a_ny_n + (1-i) \sum_{i=1}^{n-2}a_iy_i)y_n = 0 \, \, {\rm and} \, \,
2a_ny_n + (1-i) \sum_{i=1}^{n-2}a_iy_i = 0, \]
since $B$ is a domain.  From linear independence it follows that
$y_ny_{n-1} = a_{n-1}y_{n-1}y_n$, and $y_n$ is a normal element of $C,$
which is a contradiction because $y_n$ is not normal in $C$.

Case 3:  Suppose that both $g_1$ and $g_2$ are mystic
reflections that share only one of the non-invariant eigenvectors
$y_n$. Without loss of generality we may
assume that $g_1(y_n)=g_2(y_n)=-i y_n$ (by replacing $g_s$ by
$g_s^{-1}$ if necessary). Again by Definition
\ref{xxdefn1.3}(b) we may take another non-invariant eigenvector
of $g_1$ as $y_{n-1}$ such that $g_1(y_{n-1}) = i y_{n-1}$ and
the subalgebra generated by $y_{n-1}$ and $y_n$ is
isomorphic to $k\langle y_{n-1}, y_n\rangle/(y_n^2-y_{n-1}^2)$.
Since $g_1$ and $g_2$  share only the one  non-invariant eigenvector
$y_n$, $g_2(y_{n-1}) = y_{n-1}$.  In this case
$$y_{n-1}^2=g_2(y_{n-1}^2) = g_2(y_n^2)=(- iy_n)^2
= -y_n^2=-y_{n-1}^2,$$
which is a contradiction.

Case 4:  Suppose that both $g_1$ and $g_2$ are
mystic reflections that share the pair of non-invariant eigenvectors
that we can take to be $y_{n-1}$ and $y_n$. By
Definition \ref{xxdefn1.3}(b) we may assume that $g_1(y_{n-1})
=iy_{n-1}$ and $g_1(y_n)=-iy_n$. Then either $g_2(y_{n-1})=iy_{n-1}$
and $g_2(y_n)=-iy_n$, in which case $g_2=g_1$, or
$g_2(y_{n-1})=-iy_{n-1}$ and $g_2(y_n)=iy_n$, in which case $g_2=g_1^3$.
In either case we have contradicted the minimality of the set of
generators.
\end{proof}

Let $o(g)$ denote the order of $g$.

\begin{proposition}
\label{xxprop5.2}
Let $B$ be a quantum polynomial ring and let $G$ be a finite
abelian subgroup of $\Aut(B)$
generated by a minimal set of generators
$$S =\{g_1, g_2, \ldots , g_p, h_1, h_2, \ldots, h_q\}$$
where each $g_i$ is a reflection and each $h_j$ is a
mystic reflection.
\begin{enumerate}
\item
There is a $k$-linear basis
$$\{x_1,x_2, \ldots, x_p,y_1,z_1,y_2,z_2,
\ldots y_q,z_q, w_1, \ldots, w_m\}$$
of $B_1$ such that $g_i(x_i) = \lambda_i x_i$ for
$i=1,2,\ldots, p$ and $g_i(x) =x$ if $x$
is any of the other elements of the basis, and
$h_j(y_j) =iy_j, h_j(z_j)=-iz_j$ for $j=1, 2, \ldots, q$
and $h_j(x) = x$ if $x$ is any other element of the basis.
\item
The group $G$ decomposes as
$$G = \langle g_1\rangle \oplus \cdots \oplus \langle g_p\rangle \oplus
\langle h_1\rangle \oplus \cdots \oplus \langle h_q\rangle.$$
\item
The algebra $B$ is expressed as a free module over $B^G$ by
$$B = \bigoplus
x_1^{r_1}x_2^{r_2}\cdots x_p^{r_p}y_1^{s_1}z_1^{t_1}y_2^{s_2}z_2^{t_2}
\cdots y_q^{s_q}z_q^{t_q}B^G $$
where $(s_j,t_j) \in \{(0,0),(1,0),(0,1),(2,0)\}$
for $j=1, 2, \ldots , q$ and $0 \leq r_i < o(g_i)$
for $i = 1, 2, \ldots, p$.
\item
The fixed subring $B^G$ is regular.
\end{enumerate}
\end{proposition}

\begin{proof}
Since $G$ is abelian, there is a basis for $B_1$ for which the
action of each element of $G$ on $B_1$ is diagonal.  Then, by
Lemma \ref{xxlem5.1}, the elements of $S$ do not share any
non-invariant eigenvectors, and part (a) follows.

The remainder of the proof will be by induction on $p+q$.  We
will consider the generators of $G$ ordered with the reflections
listed first, followed by the mystic reflections. If $G =
\langle g_1\rangle$, then $ B= \oplus_{k=1}^{o(g_1)} x_1^kB^G$
by the proofs of \cite[Lemma 5.1]{KKZ1} for $o(g_1) >2$
and \cite[Lemma 5.2]{KKZ1} for $o(g_1)=2$. If
$G = \langle h_1\rangle,$ then $B = B^G \oplus y_1B^G \oplus
y_1^2B^G \oplus z_1B^G$ by the proof of \cite[Proposition 4.3]{KKZ1}.

Now suppose that the result holds whenever $|S| <p+q$.  First
suppose that $S=\{g_1,g_2, \ldots , g_p\}$ and let
$H = \langle g_1, g_2, \ldots, g_{p-1}\rangle$;
that is, consider the case where $q=0$. By induction we may
assume that
$$H = \langle g_1\rangle \oplus \langle g_2\rangle \oplus
\cdots \oplus \langle g_{p-1}\rangle$$ and that
$$ B= \bigoplus x_1^{r_1} x_2^{r_2}\cdots
x_{p-1}^{r_{p-1}}B^H.$$
Since $G$ is abelian, $g_p(B^H)=B^H$. By part (a) we have that
$h(x_p) = x_p$ for all $h \in H,$ and thus $x_p \in B^H.$
Since $g_p(x_p) = \lambda_p x_p$ where $\lambda_p$ is an $o(g_p)$th
root of unity, $g_p^k(x_p)\neq x_p$ for $1 \leq k < o(g_p)$.
Consequently, $g_p^k \notin H$ for $1 \leq k < o(g_p),$ and
$G = \langle g_1\rangle \oplus \cdots \oplus \langle g_p\rangle.$
Clearly $B^G = B^H\cap B^{g_p}$.  As above
$$B = \oplus_{k=1}^{o(g_p)} x_p^kB^{g_p}.$$
Thus if $a \in B^H,$ then $a = \sum_{k=1}^{o(g_p)} x_p^k a_k$
for unique $a_k$'s in $B^{g_p}$.  Let $h \in H$.  Then
$$ h(a) = a = \sum_{k=1}^{o(g_p)} h(x_p^k) h(a_k) =
\sum_{k=1}^{o(g_p)} x_p^k h(a_k).$$
Since $G$ is abelian, we have that $h(a_k) \in B^{g_p},$
and  from the uniqueness of the sum, that $h(a_k)=a_k.$
Consequently each $a_k \in B^{g_p}\cap B^H =
B^G$, and $ B^H =\oplus_{k=1}^{o(g_p)} x_p^kB^G.$
It follows that
$$ B= \bigoplus x_1^{r_1} x_2^{r_2}\cdots
x_{p-1}^{r_{p-1}}B^H = \bigoplus x_1^{r_1} x_2^{r_2}\cdots
x_{p-1}^{r_{p-1}}x_p^{r_p}B^G.$$

Now suppose that $S =\{g_1, g_2, \ldots , g_p, h_1, h_2, \ldots, h_q\}$
with $q \geq 1,$ and let
$$H = \langle g_1, g_2, \ldots , g_p, h_1, h_2, \ldots, h_{q-1}\rangle.$$
Then by induction we have that
$$H = \langle g_1\rangle \oplus \cdots \oplus \langle g_p\rangle \oplus
\langle h_1\rangle \oplus \cdots \oplus \langle h_{q-1}\rangle,$$
and $B$ is a free module over $B^H$ written as
$$B = \bigoplus x_1^{r_1}x_2^{r_2}\cdots
x_p^{r_p}y_1^{s_1}z_1^{t_1}y_2^{s_2}z_2^{t_2}
\cdots y_{q-1}^{s_{q-1}}z_{q-1}^{t_{q-1}}B^H. $$
By part (a) we have that $y_q \in B^H.$  The order of $h_q$
is $4$, and since $h_q(y_q)=iy_q, h_q^2(y_q)=-y_q$, and
$h_q^3(y_q)=-iy_q$, we have that no nonidentity power of
$h_q$ is in $H$. Thus $\langle h_q\rangle \cap H = \{e\}$ and
$$G = \langle g_1\rangle \oplus \cdots \oplus \langle g_p\rangle \oplus
\langle h_1\rangle \oplus \cdots \oplus \langle h_q\rangle.$$
As noted above $B$ as a free module over $B^{h_q}$ decomposes as
$$B= B^{h_q}\oplus y_qB^{h_q}\oplus y_q^2B^{h_q}\oplus z_qB^{h_q}.$$
Using that $B^G = B^{h_q}\cap B^H,$ and the inductive decomposition
of $B$ as a free module over $B^H$ yields that
$$B = \bigoplus x_1^{r_1}x_2^{r_2}\cdots
x_p^{r_p}y_1^{s_1}z_1^{t_1}y_2^{s_2}z_2^{t_2}
\cdots y_q^{s_q}z_q^{t_q}B^G.$$
Thus parts (b,c) follow by induction.

(d) This follows from \cite[Lemma 1.10]{KKZ1}.
\end{proof}

\begin{theorem}
\label{xxthm5.3}
Let $B$ be a quantum polynomial ring and let $G$ be a finite
abelian subgroup of $\Aut(B)$.  Then $B^G$ is regular
if and only if $G$ is generated by quasi-reflections.
\end{theorem}

\begin{proof} One implication is Proposition \ref{xxprop5.2}(d)
and the other implication is Proposition \ref{xxprop1.5}(a).
\end{proof}

By Proposition \ref{xxprop5.2}(b) Conjecture \ref{xxcon4.8} holds
in this case.

In the rest of this section we give an example
that indicates that we might need to go to the world
of Hopf algebras for the most general version of the Shephard-Todd-Chevalley
Theorem.

\begin{example}
\label{xxex5.4}
Let $C$ be the quantum $2\times 2$-matrix algebra
${\mathcal O}_q(M_2)$ \cite[Definition I.1.7]{BG}
that is generated by $x_{11},x_{12},x_{21},x_{22}$
subject to the relations
$$
\begin{aligned}
x_{12}x_{11}&=qx_{11}x_{12}\\
x_{21}x_{11}&=qx_{11}x_{21}\\
x_{22}x_{12}&=qx_{12}x_{22}\\
x_{22}x_{21}&=qx_{21}x_{22}\\
x_{21}x_{12}&=x_{12}x_{21}\\
x_{22}x_{11}&=x_{11}x_{22}+(q^{-1}-q)x_{12}x_{21}.
\end{aligned}
$$
This algebra is in fact a bialgebra, with coproduct
$\Delta$ and counit $\epsilon$  given by:
\begin{align}
\label{5.4.1}
\Delta(x_{ij})&=\sum_{s=1}^2 x_{is}\otimes x_{sj}
\tag{5.4.1}
\\
\epsilon(x_{ij})&=\begin{cases}
1& i=j\\ 0&i\neq j\end{cases}
\notag
\end{align}
for all $i,j\in [2]$.
\end{example}

\begin{lemma}
\label{xxlem5.5}
Let $C={\mathcal O}_q(M_2)$ be defined as in Example \ref{xxex5.4}.
\begin{enumerate}
\item
The algebra $C$ is an iterated Ore extension
$k[x_{11}][x_{12};\tau_2][x_{21};\tau_3][x_{22};
\tau_4,\delta_4]$, and it is a connected graded
algebra with $\deg x_{ij}=1$.
\item
It is a quantum polynomial ring of
dimension 4.
\item
It is an Auslander regular, Cohen-Macaulay,
noetherian domain.
\item
Suppose $q\neq \pm 1$. Then every normal element in degree 1
is of the form $c_1x_{12}+c_2x_{21}$ for $c_i \in k$.
As a consequence, $C$ is not isomorphic to a skew
polynomial ring.
\end{enumerate}
\end{lemma}

\begin{proof} (a) \cite[Example I.1.16]{BG}.

(b,c) These follows from part (a) and \cite[Lemma I.15.4]{BG}.

(d) This can be checked directly using
the relations listed above.
\end{proof}

In the next lemma and the following proposition let $q$ be a
primitive $m$th root of unity where $m\geq 3$.
Let
$$n=\begin{cases} m& {\text{if $m$ is odd}}
\\ \frac{m}{2} & {\text{if $m$ is even}.}
\end{cases}$$
Let $H$ be a factor bialgebra of $C$ modulo $x_{11}^n-1, x_{12}^n,
x_{21}^n, x_{22}^n-1$. Let $C(n)$ be the subspace of $C$ spanned by
$x_{11}^ax_{12}^bx_{21}^cx_{22}^d$ for all $0\leq a,b,c,d\leq n-1$.

\begin{lemma}
\label{xxlem5.6}
Let $H=C/I$ where
$I$ is the ideal generated by $x_{11}^n-1, x_{12}^n, x_{21}^n,
x_{22}^n-1$.
\begin{enumerate}
\item
Then $H$ is a finite dimensional Hopf algebra
with the coalgebra structure determined by \eqref{5.4.1}.
\item
The canonical map $\phi: C\to H$ induces an isomorphism
of vector spaces $C(n)\cong H$.
\item
Let $y_{ij}$ be the image of $x_{ij}$ in $H$. Then $H$ has
a $k$-basis
$$\{y_{11}^ay_{12}^by_{21}^cy_{22}^d\;|\;
0\leq a,b,c,d\leq n-1\}.$$
\end{enumerate}
\end{lemma}

\begin{proof} (a) First we prove that $H$ is a bialgebra.
It suffices to show that $I$ is a bialgebra ideal of
$C$. Note that
$$\Delta(x_{ij})=
x_{i1}\otimes x_{1j}+x_{i2}\otimes x_{2j},$$
and
$$(x_{i2}\otimes x_{2j})(x_{i1}\otimes x_{1j})
=q^2(x_{i1}\otimes x_{1j})(x_{i2}\otimes x_{2j}).$$
Since $q^2$ is a primitive $n$th root of unity,
$$\Delta(x_{ij}^n)=x_{i1}^n\otimes x_{1j}^n+x_{i2}^n
\otimes x_{2j}^n.$$
Hence
\begin{align}
\label{5.6.1}
\Delta(x_{11}^n-1)&=(x_{11}^n-1)\otimes x_{11}^n+1\otimes (x_{11}^n-1)
+x_{12}^n \otimes x_{21}^n,
\tag{5.6.1}\\
\Delta(x_{12}^n)&=x_{11}^n\otimes x_{12}^n+x_{12}^n
\otimes x_{22}^n,\notag\\
\Delta(x_{21}^n)&=x_{21}^n\otimes x_{11}^n+x_{22}^n
\otimes x_{21}^n,\notag\\
\Delta(x_{22}^n-1)&=(x_{22}^n-1)\otimes x_{22}^n+1\otimes (x_{22}^n-1)
+x_{21}^n \otimes x_{12}^n.\notag
\end{align}
Therefore $I$ is a bialgebra ideal and $H$ is a bialgebra.
It remains to show that $H$ is a Hopf algebra.
Let $B$ be the localization $C[D_q^{-1}]$, where
$D_q$ is the quantum determinant $x_{11}x_{22}-q^{-1}x_{12}x_{21}$
(which is a central element in $C$). Then $B$ is
a Hopf algebra denoted by ${\mathcal O}_q(GL_2)$
\cite[Definition I.1.10]{BG}. Since $D_q^{m}$ is
equal to $1$ in $H$, $H$ is a factor bialgebra of $B$
modulo $J=(x_{11}^n-1, x_{12}^n, x_{21}^n,x_{22}^n-1)$.
Using the definition of the antipode \cite[Definition I.1.10(7)]{BG},
it is easy to see that $J$ is a Hopf ideal of $B$.
Thus $H$ is a Hopf algebra.

(b,c) Clear.
\end{proof}

We define a right $H$-comodule structure on $C$ by
$$\rho(x_{ij})= \sum_{s=1}^2 x_{is}\otimes y_{sj}$$
for all $i,j\in [2]$.

By the next proposition, we may think of $H$ as a quantum
reflection group of $C$. By the definition of $n$, $q^{n^2}$
is either 1 or $-1$.

\begin{proposition}
\label{xxprop5.7}
The algebra $C$ defined above is a right $H$-comodule algebra,
and the coinvariant subring $C^{co H}$ is the skew
polynomial ring $k_{\pm 1}[x_{11}^n,x_{12}^n,x_{21}^n,x_{22}^n]$
subject to the relations
$$
\begin{aligned}
x_{12}^n x_{11}^n &=q^{n^2}x_{11}^n x_{12}^n \\
x_{21}^n x_{11}^n &=q^{n^2}x_{11}^n x_{21}^n \\
x_{22}^n x_{12}^n &=q^{n^2}x_{12}^n x_{22}^n \\
x_{22}^n x_{21}^n &=q^{n^2}x_{21}^n x_{22}^n \\
x_{21}^n x_{12}^n &=x_{12}^n x_{21}^n\\
x_{22}^n x_{11}^n &=x_{11}^n x_{22}^n.
\end{aligned}
$$
\end{proposition}

\begin{proof}
Since $C$ is a bialgebra, it is a right $C$-comodule
algebra. Since $H$ is a factor bialgebra of $C$, $C$ is
also a right $H$-comodule algebra.

Let $B$ be the subring of $C$ generated by $x_{11}^n,x_{12}^n,
x_{21}^n,x_{22}^n$. Then $B$ is isomorphic to
${\mathcal O}_{q^{n^2}}(M_2)$, which is a skew polynomial ring
$k_{\pm 1} [x_{11}^n,x_{12}^n,x_{21}^n,x_{22}^n]$ described in
the statement. Using \eqref{5.6.1} we see that
$$\begin{aligned}
\rho(x_{11}^n-1)&=(x_{11}^n-1)\otimes 1,\\
\rho(x_{12}^n)&=x_{12}^n
\otimes 1,\\
\rho(x_{21}^n)&=x_{21}^n\otimes 1,\\
\rho(x_{22}^n-1)&=(x_{22}^n-1)\otimes 1.
\end{aligned}
$$
This implies that $B\subset C^{co H}$. It remains to show that
$C^{co H}\subset B$. Let $f$ be a nonzero element in $C^{co H}$. Write
$f=\sum h_i g_i$, where $h_i$ are linear combinations
of monomials of the form $x_{11}^ax_{12}^bx_{21}^cx_{22}^d$
for $0\leq a,b,c,d<n$, and where $\{g_i\in B\}$ are
linearly independent over $k$. Then
$$\sum_i h_ig_i\otimes 1=f\otimes 1=\rho(f)=\sum_i \rho(h_i)\rho(g_i)=\sum_i
\rho(h_i) (g_i\otimes 1).$$
Hence $\sum_i (h_i\otimes 1-\rho(h_i))(g_i\otimes 1)=0$,
where $g_i\otimes 1\in B\otimes k$ and $h_i\otimes 1
-\rho(h_i)\in C(n)\otimes H$. Recall that $C(n)=\sum_{0\leq a, b, c, d<n}
kx_{11}^ax_{12}^bx_{21}^cx_{22}^d$. Write
$h_i\otimes 1-\rho(h_i)=\sum_{0\leq a, b, c, d<n}
x_{11}^ax_{12}^bx_{21}^cx_{22}^d\otimes w^i_{a,b,c,d}$.
Then $\sum_{0\leq a,b,c,d<n,i} x_{11}^ax_{12}^bx_{21}^cx_{22}^d g_i
\otimes w^i_{a,b,c,d}=0$ implies that $w^i_{a,b,c,d}=0$ for all
$a,b,c,d,i$ because $\{x_{11}^ax_{12}^bx_{21}^cx_{22}^d g_i\}$ are
linearly independent over $k$.
This means that $\rho(h_i)=h_i\otimes 1$ or $h_i\in C^{co H}$ for all $i$.
Since $H$ is a factor bialgebra of $C$, $\rho(h_i)=h_i\otimes 1$
implies that $\Delta_H(\phi(h_i))=\phi(h_i)\otimes 1$
where $\phi:C\to H$ is the canonical map. Hence $\phi(h_i)\in k$.
By Lemma \ref{xxlem5.6}(b), $\phi: C(n)\to H$ is an isomorphism of vector
spaces which sends $x_{11}^ax_{12}^bx_{21}^cx_{22}^d$ to
$y_{11}^ay_{12}^by_{21}^cy_{22}^d$. By Lemma \ref{xxlem5.6}(c),
$h_i\in k$. Therefore $C^{co H}=B$.
\end{proof}

The next proposition is a Shephard-Todd-Chevalley Theorem for the algebra $C$.
It also shows that $C$ is a rigid graded algebra, i.e. $C^G \not\simeq C$
for all finite groups $G$ of graded automorphisms of $C$. Combined
with Proposition \ref{xxprop5.7} it shows that the
``quantum reflection group'' $H$ of $C$ provides a regular ring of
invariants that is different from the fixed subrings under any group
action. This suggests that Hopf actions may provide the proper
context for a generalized Shephard-Todd-Chevalley Theorem.

\begin{proposition}
\label{xxprop5.8}
Let $C$ be as defined above. Suppose $q\neq \pm1 $.
\begin{enumerate}
\item
$C$ has no mystic reflections.
\item
All quasi-reflections of $C$ are reflections of the form
$$g_b(x_{11}) = x_{11}, \; \; g_b(x_{12}) = b x_{21}, \;
\; g_b(x_{21}) = b^{-1}x_{12}, \;\; g_b(x_{22}) = x_{22}$$
for $b \in k^{\times}$.
\item
The finite groups $G\subset \Aut(C)$ that are generated by
quasi-reflections of $C$ are the dihedral groups and the
cyclic group of order 2.
\item
For all finite groups $G\subset \Aut(C)$ that are generated
by quasi-reflections of $C$, the fixed subring $C^G$ is regular
and is generated as an algebra by 3 elements.  Hence
$C^G \not\cong C$ for any non-trivial finite group $G\subset \Aut(C)$.
\item
Conjectures \ref{xxcon0.2} and \ref{xxcon4.8} hold for $C$.
\item
Assume that $q$ is a root of unity. Then for any finite group $G\subset \Aut(C)$,
$C^G \not\cong C^{co H}$ for the Hopf algebra $H$ in
Proposition \ref{xxprop5.7}.
\end{enumerate}
\end{proposition}

\begin{proof}
(a) If $\sigma$ is a mystic reflection, then there are
$z_1,z_2\in C_1$ such that $z_1z_2=-z_2z_1$. But
it is easy to check that this is impossible when
$q\neq -1$.

(b) By Lemma \ref{xxlem5.5}(d) every normal element of degree 1
is in $kx_{12}+kx_{21}$.

Let $g$ be a quasi-reflection that is a reflection.  Then there
is a normal element $y \in C_1$ with
$g(y) = \lambda y\in kx_{12}+kx_{21}$ and
a basis $\{y_1,y_2,y_3,y_4\}$ of $C_1$ with $y_1=y$ and $g(y_i) = y_i$
for $i\neq 1$,  Then $g$ induces a graded automorphism that is the
identity map on the factor algebra $C/\langle y\rangle$, and so
$$g(x_{12}) = x_{12} + \alpha y, \;\;
  g(x_{21}) = x_{21} + \beta y,\;\;
  g(x_{11}) = x_{11} + \gamma y,\;\;
  g(x_{22}) = x_{22} + \delta y $$
for scalars $\alpha, \beta, \gamma, \delta \in k$.  Since
$g(x_{12}x_{11}) = q g(x_{11}x_{12})$ we get either $\gamma = 0$
or $\alpha y = -x_{12}$.  Suppose that $\gamma \neq 0$.
Then $\alpha \neq 0$. In a similar way, we have that
$g(x_{21}x_{11}) = q g(x_{11}x_{21})$ yields $\gamma = 0$ or
$\beta y = -x_{21}$.  If $\gamma \neq 0$ then $\beta \neq 0$, and
$y = -x_{12}/\alpha = -x_{21}/\beta$ is a contradiction.  Hence
$\gamma = 0$. By symmetry $\delta = 0$, and whence $x_{11}$ and
$x_{22}$ are fixed by $g$.  The last quadratic relation implies
that $g(x_{12}x_{21}) = x_{12}x_{21}$.
To show that $g$ is of the form indicated, write
$y = ax_{12} + bx_{21}$ and compute
$g(x_{12}x_{21}) = x_{12}x_{21}
= (x_{12} + \alpha y) (x_{21} + \beta y)$.  Equating terms
gives the result.

(c) If $G$ is generated by one reflection, then $G$ is cyclic of
order 2.  Suppose $G$ has at least two generators $g_a$ and $g_b$.
Then $g_a g_b = d_{{ab}^{-1}}$, where $d_c$ is the ``diagonal map''
determined by
$$d_c(x_{11} )= x_{11}, d_c(x_{12})= cx_{12},
d_c(x_{21})= c^{-1}x_{21}, d_c(x_{22} )= x_{22}.$$
Let $N=\{d_c \in G\}$; then $N$ is a normal subgroup of $G$ with
$[G:N]=2$.  The group $N$ is isomorphic to a finite subgroup of
$k^{\times}$, so is cyclic, generated by an $m$th root of unity
$\lambda$.  Furthermore, $g_1$ and $N$ generate $G$, and $G$ is
isomorphic to the dihedral group $D_m$ of order $2m$.

(d,e,f) It is easily checked that the fixed subring under the
dihedral group $D_m$ represented as above is generated as an
algebra by $x_{11},  x_{22},  x_{12}x_{21}, $ and  $x_{12}^m+b^mx_{21}^m$.
When $G=\langle g_b\rangle $, it is generated by $x_{11}, x_{22},
x_{12}x_{21}, $  and  $x_{12}+bx_{21}$. Since
$$x_{12}x_{21} =(q^{-1}-q)^{-1}(x_{22} x_{11}-x_{11}x_{22})$$
the fixed subring $C^G$ is generated by 3 elements in both cases.
It can be checked that it is an iterated Ore-extension, so is regular.
As both $C$ and $C^{co H}$ require four algebra generators, $C^G$ is
not isomorphic to either $C$ or $C^{co H}$.

By Proposition \ref{xxprop1.5}(b) $C^G$ is regular (or has finite
global dimension) if and only if $G$ is generated by quasi-reflections,
so Conjecture \ref{xxcon0.2} holds. Conjecture \ref{xxcon4.8} holds
trivially since $D_{2m}$ is a reflection group in the classical
sense. Thus we have proven part (e).

Finally if $C^G$ is isomorphic to either $C$ or $C^{co H}$ for
some finite group $G\subset \Aut(C)$, then
$C^G$ is regular. By part (e) $G$ is either trivial or generated
by quasi-reflections. But $C^G\not\cong C$ or $ C^{co H}$ if $G$ is generated
by quasi-reflections, as proved in the first paragraph of the
proof of (d,e). Hence $G$ is trivial. In particular, $C$ is rigid.
\end{proof}

Proposition \ref{xxprop0.5} follows from Propositions
\ref{xxprop5.7} and \ref{xxprop5.8} immediately. Note
that a Hopf algebra $H$ coaction  on $C$ is equivalent to a
Hopf algebra $H^*$ action on $C$. Neither the Hopf algebra $H$
in Lemma \ref{xxlem5.6}, nor its dual $H^*$, is semisimple.

\section{The groups $M(n,\alpha,\beta)$}
\label{xxsec6}

Classical (pseudo-)reflection groups are well-understood in many
respects. As generalizations of reflection groups, the groups $M(n,\alpha,\beta)$,
generated by mystic reflections,
merit further investigation. We begin their study here.

In this section we give some examples of groups  $M(n,\alpha,\beta)$
which answer Question \ref{xxque0.4}(a).
For simplicity let $k={\mathbb C}$.
In the first example, we give some details about groups
$M(n,1,\beta)$ generated by mystic reflections. These groups may
be compared to the classical reflection groups classified by
Shephard and Todd \cite{ShT}.  In the classical case there are
three infinite families of reflection groups (i.e. groups
with a representation generated by reflections of
$k[x_1, \cdots, x_n]$ for some $n$):
\begin{equation}
\label{xxequ6.0.1}
\begin{cases}
{\text{the cyclic groups}}, & \\
{\text{the symmetric groups}}, & {\text{and}} \\
{\text{the groups $G(m,p,n)$}} &{\text{(which include the dihedral groups)}};
\end{cases}
\tag{6.0.1}
\end{equation}
there are also 34 exceptional groups.
The group $G(m,p,n)$ (for positive integers $m,p,n$, where $p$
divides $m$, so $m=pq$) is a semidirect product
of the set of $n \times n$ diagonal matrices $A(m,p,n)$ described
below, and $S_n$ represented as permutation matrices. The group
$A(m,p,n)$ is
\[A(m,p,n) =
\left\{
\matfive{ \omega_1 & 0 & \cdots & 0 & 0\\
                0 &\omega_2 & \cdots & 0 & 0\\
                \vdots & \vdots & \ddots & \vdots & \vdots\\
                0 & 0 & \cdots & \omega_{n-1} & 0\\
                0 & 0 & \cdots & 0 & \omega_n}: w_i^m = 1
                \text{ and } (\omega_1 \cdots \omega_n)^{m/p} = 1
  \right\}.
  \]
Then  $G(m,p,n)$ is a subgroup of weighted $n \times n$
permutation matrices, and $G(m,p,n)$ contains a copy of the
$n \times n$ permutation matrices.  The order of $G(m,p,n)$ is
$m^nn!/p$ (see e.g. \cite[pp. 161-2 and p.166]{K} and \cite{ShT}).
The group $G(m,m,2)$ is isomorphic to the dihedral group with
$2m$ elements, and for arbitrary $p$ the groups $G(m,p,2)$ are
realized as symmetry groups of certain complex polytopes.

\begin{example}
\label{xxex6.1}
Let $A$ be the skew-polynomial ring $\mathbb{C}_{-1}[x_1,x_2,x_3]$
with $p_{ij} = -1$ for all $i \neq j$. Let $G$ be the group
generated by the two mystic reflections $\tau_{1,2,1}$ and
$\tau_{2,3,1}$; then $G = M(3,1,2)$.  We first note that $G$ is
the rotation group of cube,  viewing a cube centered at the origin,
with the three coordinate axes through the center of the faces.
Then the rotation group of the cube is isomorphic to the
symmetric group $S_4$ and is generated by the two matrices:
\[g_1:= \matthree{ 0 & -1 & 0\\ 1 & 0 & 0\\ 0 & 0 & 1} \;\;\hbox{ and }
\; g_2 := \matthree{ 1 & 0 & 0\\ 0 & 0 & -1\\ 0 & 1 & 0}\]
that act on $A$ as the mystic reflections $g_1 = \tau_{1,2,1}$
and $g_2 = \tau_{2,3,1}$, respectively.   By Proposition
\ref{xxprop4.4}(d) the ring of invariants $A^G$ is the
commutative polynomial ring
\[\mathbb{C}[x_1^2+x_2^2+x_3^2,\; x_1x_2x_3,\; x_1^4+x_2^4+x_3^4].\]
Hence $S_4$ is a ``reflection group" of $A$ in the sense that it is
has a representation that is generated by quasi-reflections of $A$.
However, as in Proposition \ref{xxprop4.4}(c), this representation
of $S_4$ contains no non-identity reflections of the commutative
polynomial ring $\mathbb{C}[x_1,x_2,x_3]$ under the usual action.
It is interesting that the generators of $A^G$ in the
case of quantum polynomials are nicer than those for the
commutative polynomial ring $\mathbb{C}[x_1,x_2,x_3]^G$,
see \cite[p. 337, Problem 12]{CLO}.
Although the representation of $S_4$ generated by the mystic
reflections above is not a ``reflection group" in the commutative
sense (the group of matrices above is not generated by classical
reflections) the group $S_4$ has another representation (the
representation as permutation matrices) that is a
``reflection group" in the classical sense.
\end{example}

Our next examples show that there exist infinite families of
finite groups $G$ having a representation that is generated by
quasi-reflections of a regular algebra, yet the abstract group $G$ has
\underline{no} representation generated by classical reflections.
Though not classical reflection groups, these groups are
``reflection groups" in the sense that they have representations
that act on a regular algebra  producing a regular fixed ring
(the fixed ring is even a commutative polynomial ring).

\begin{example}
\label{xxex6.2}
We claim that the mystic reflection groups $M(2,1,2m)$, for $m\gg 0$, are
not isomorphic to classical reflection groups as abstract groups.

Let $A$ be the quantum plane ${\mathbb C}_{-1}[x_1,x_2]$ and
let $G$ be the group generated by the mystic reflections
$g_1= \tau_{2,1,1}$ and $g_2=\tau_{2,1,\lambda}$ where
$\lambda$ is a primitive $2m$th root of unity.  Then
\[g_1= \mattwo{ 0 & 1\\ -1 & 0}\;\; \text{ and }\;\;
g_2 = \mattwo{ 0 & \lambda\\ -\lambda^{-1} & 0}.\]
The group $G$ is also generated by $g_1$ and
\[g_3 = \mattwo{ \lambda & 0\\ 0 & \lambda^{-1}},\]
and the fixed ring is $A^G =
\mathbb{C}[x_1x_2, x_1^{2m}+x_2^{2m}]$.
Then taking $a = g_3$ and $b = g_1$, it is not difficult
(e.g. \cite[Exercise 16, p. 25-26]{N}) to see that
$G = M(2,1, 2m)$ is isomorphic to the ``dicyclic group"
$Q_{4m}$ generated by $a$ and $b$ with relations:
\[ a^{2m} = 1, \hspace*{.2in}  b^{-1} a b = a^{-1},
\hspace*{.2in}  b^2 = a^m.\]
When $m=1$ then $Q_4$ is the cyclic group of order $4$, and
when $m=2$ then $Q_8$ is the quaternion group of order $8$.
It is easy to see that all elements
of $G$ can be written in the form $a^i$ or $a^ib$, that $|G| = 4m$, and
that $Q_{4m}$ has a unique element $a^m$ of order $2$.  When $m=2^{n-2}$
then the dicyclic groups are called ``generalized quaternion groups"
$Q_{2^n}$, generated by elements $a$ and $b$ with relations:
\[a^{2^{n-1}} = 1, \hspace*{.2in} b^{-1}ab = a^{-1},
\hspace*{.2in} b^2 = a^{2^{n-2}}\]
for $n \geq 3$ (e.g. \cite[Exercise 15, p. 25]{N}).

As each of the 34 exceptional complex reflection groups has order divisible
by $4$, it is possible that a finite number of  $M(2,1,2m) \cong Q_{4m}$ could be isomorphic to
an exceptional reflection group; however, by examining the list of the 34 exceptional
groups and considering their orders, one can easily determine that e.g. $Q_8, Q_{12}, Q_{16}, Q_{20}, Q_{28}$,
or any of the generalized quaternion groups cannot be isomorphic to an exceptional reflection group.

We claim that except when $m = 1$ (when $G \cong \mathbb{Z}/4\mathbb{Z}$),
the groups $Q_{4m}$ are not isomorphic to any of the groups in
the three infinite families of groups in the Shephard-Todd
table of complex reflection groups \eqref{xxequ6.0.1}.
These non-abelian groups $Q_{4m}$ each
have a unique element of order $2$, so are not isomorphic to cyclic,
symmetric or dihedral groups.    Hence,
once we show that $Q_{4m}$ is never
isomorphic to a group of the form $G(m,p,n)$, there are
an infinite number of new ``reflection groups" arising from reflections of
${\mathbb C}_{-1}[x_1,x_2]$.

Finally, we show that the groups $Q_{4m^*}$ (we changed $m$ to a
different integer $m^*$) are not isomorphic to a group $G(m,p,n)$.
Since the dicyclic groups each have a unique element of order $2$,
then if a dicyclic group $Q_{4m^*}$ is
isomorphic to the reflection group $G(m,p,n)$, it follows that
$n = 2$ because it contains a subgroup isomorphic to $S_n$.
The groups $G(m,p,2)$ are generated by the transposition
\[t = \mattwo{0 & 1\\ 1 & 0}\]
and the matrices
\[\mattwo{\omega^{a_1} & 0 \\ 0 & \omega^{a_2}}\]
with $a_1 + a_2 \equiv 0 \mod p$, where $\omega$ is a primitive
$m$th root of unity \cite[p. 147]{N}.  In order that
$|G(m,p,2)| = (m^2)2!/p = 4 m^* = |Q_{4m^*}|$, $m$ must be even,
and hence the element
\[\mattwo{-1&0\\0&-1} =
\mattwo{\omega^{m/2} & 0 \\ 0 & \omega^{m/2}},\]
as well as the transposition $t$, are two elements of $G$ with
order $2$, contradicting the uniqueness of the elements of order
$2$ in $Q_{4m^*}$.  Hence it follows that the groups $G= M(2,1,2m)$
are not isomorphic to the groups $G(m,p,n)$, and so can be
classical reflection groups only if they happen to be one of the
34 exceptional complex classical reflection groups in the
Shephard-Todd table.

Therefore we proved the assertion by taking
$m^* \gg 0$.
\end{example}

In Example \ref{xxex6.1} we showed that $M(3,1,2)$ is isomorphic
to the group of rotations of the cube.  Next we examine the family of groups
$G=M(n,1,2)$ that are generated by mystic reflections of
${\mathbb C}_{-1}[x_1, \cdots, x_n]$.  It is easy to check that
$G$ has some properties of the classical reflection group
$G'= G(2,2,n)$; both groups have order $2^{n-1} n!$ and
both contain the diagonal matrices $D = A(2,2,n)$, described
above, as a normal subgroup. In the rest of this section we show
that when $n$ is even the groups $G$ and $G'$ are not isomorphic,
but when $n$ is odd they are isomorphic.  Hence $M(n,1,2)$ for
$n$ even  provides a second family of new ``reflection groups"
generated by mystic reflections.

\begin{example}
\label{xxex6.3}
The group $G= M(n,1,2)$ is generated by the mystic reflections
$\tau_{i,j,1}$ for all $i \neq j$.  Both $G$ and $G'$ are
subgroups of the group $B$ of all signed ($\pm1$) permutation
$n \times n$ matrices.  The order of $B$ is $2^n n!$, and both
$G$ and $G'$ have order $2^{n-1}n!$.  The group $G$ is the
kernel of the determinant map from $B \rightarrow \{\pm1\}$.
The diagonal matrices in $G$ are precisely the diagonal matrices
$D= A(2,2,n)$, described above, as a subgroup of $G(2,2,n)$; these
matrices contain an even number of entries of $-1$ and are the
diagonal matrices  in $B$ with determinant $1$.  The subgroup $D$
is normal in $B$ (and hence in both $G$ and $G'$) and $G/D$ is
isomorphic to $S_n$ by reducing $-1$ entries to $1$ entries;
it is clear from its definition that $G'/D$ is isomorphic to
$S_n$. Hence, it is not surprising that some effort is required
in distinguishing the groups $G$ and $G'$ when they are not 
isomorphic.\\

\noindent
Case 1: $n$ is odd.  Consider the elements $\sigma_i$,
$i = 1, \ldots, n-1$ chosen so that the coset $D\sigma_i$
represents the transposition $(i, i+1)$ in the factor
group $G/D$; since $n$ is odd we can choose $\sigma_i$ to be of
the form
 \[\sigma_i =
 \mateight{ -1 &  & & & & & & \\
                                    &  \ddots & & & & & & \\
                                     &  & -1 &  &  &  & & \\
                  & &  & 0 & 1 &  & & \\
                  & & &1 & 0 & & & \\
               &  & & & & -1& & \\
                 & & & & & &\ddots & \\
               & & & & & & & -1 }\]
               as the determinant of $\sigma_i$ is $1$.  One can
 check that the elements $\sigma_i$ of $G$ satisfy the
 relations of $S_n$:
\[\sigma_i \sigma_{i+1} \sigma_i = \sigma_{i+1} \sigma_i \sigma_{i+1}
\; \;\text{ for } i = 1, \ldots, n-1,\]
 \[\sigma_i \sigma_j = \sigma_j \sigma_i \;\; \text{  for  } |i-j|
\geq 2, \text{ and }\]
 \[\sigma_i^2 = 1 \text{ for } i = 1, \ldots, n-1.\]
 Hence $G$ contains a copy of $S_n$ that acts on $D$ as it does in $G'$,
so $G \cong G'$.\\

\noindent
Case 2: $n$ is even.  We will show that the groups $G$
and $G'$ differ in
the distribution of elements of order $2 ^{t}$ for some $t$.
First we prove the following lemmas:

\begin{lemma}
\label{xxlem6.4}
Let $M$ be a $t \times t$ weighted $t$-cycle matrix of the
 form
 \[ M = \matfive{ 0 & &\cdots &  0 &\pm 1\\
                             \pm 1 & 0 &  \cdots& &0\\
                             0 & \pm 1 & \ddots &\vdots &\vdots\\
                           \vdots  & \ddots& \ddots& 0 &\\
                            0 & \cdots &0& \pm 1 & 0}.\]
Then the order of $M$ is $t$ when the number of $-1$'s in $M$  is
an even number, and it is $2t$ when the number of $-1$'s in $M$
is an odd number.
\end{lemma}

\begin{proof}  Thinking of $M$ as a function makes this easy to check.
\end{proof}

\begin{lemma}
\label{xxlem6.5}
 Let $g$ be an element of $G'$ with $Dg$ an even permutation of
$G'/D \cong S_n$.  Then $g \in G$.
\end{lemma}

\begin{proof}
 The coset $Dg \in G'/D$ contains a permutation of determinant $1$, and
hence all elements of $Dg$
 have determinant $1$, and so are in $G$.
\end{proof}

We next compare the orders of $g$ and $Dg$ when $g$ has order a
power of $2$. Note that it is possible for the order of $Dg$ to
be half the order of $g$; e.g. if
\[g = \matfour{ 0 & -1& 0 & 0\\1 & 0 & 0 & 0\\ 0 & 0 & -1 &0
\\ 0 & 0 & 0 & -1}\]
then $g$ has determinant $1$ so is in $G$ and
$|g| = 4$, but the order of $Dg$ is $2$ because changing the
sign of the first row and the third row gives the element
\[ \matfour{ 0 & 1& 0 & 0\\1 & 0 & 0 & 0\\ 0 & 0 & 1 &0
\\ 0 & 0 & 0 & -1}\]
in the coset $Dg$ that clearly has order $2$; hence it is clear
that the coset $Dg$ has order $2$.  The next lemma shows
that the order of $Dg$ cannot decrease from $|g|$ by more than
one power of $2$.

\begin{lemma}
\label{xxlem6.6}
If the order of $g \in G$ is $2^{\ell+1}$, then the order of
the coset $Dg$ in $G/D$ is either $2^{\ell}$ or $2^{\ell+1}$.
\end{lemma}

\begin{proof}
The order of the coset $Dg$ in $G/D \cong S_n$ must be
$\leq 2^{\ell +1}$, so it can be represented by a  block
permutation matrix for a product of disjoint $2^{t_i}$-cycles
for $t_i \leq \ell +1$.  Hence $g$ is obtained from the block
permutation matrix representing the product of these
$2^{t_i}$-cycles for $t_i \leq 2^{\ell +1}$  by an even number
of sign changes.   It follows from Lemma \ref{xxlem6.4} that
if all the cycles in $Dg$ are of length $\leq 2^{\ell -1}$, then $g$
cannot have order $2^{\ell + 1}$.
\end{proof}

\begin{lemma}
\label{xxlem6.7}
Let $Dg$ be an element of $G/D\cong S_n$ of order $2^{k}$ for $k \geq \ell$ where
$n=2^\ell m$ for $m$ odd, and let
$\sigma$ be the corresponding permutation in $S_n$. Let $a$ in $G'$
be a representative of the coset corresponding to $\sigma$ in
$G'/D \cong S_n$.
Then the number of elements in $Dg$ of order $2^{k + 1}$ in $G$ is
$\geq $ the number of elements in $Da$ of order $2^{k +1}$ in $G'$.
\end{lemma}
\begin{proof}  Since $D$ is normal in $B$, and $G$ and $G'$ are
normal in $B$, conjugation by
any permutation matrix preserves the number of elements of each
order in $Dg$ and $Da$.
Hence we may assume that elements of $G$ and $G'$ have block form
corresponding
to disjoint weighted cycles of the form in Lemma \ref{xxlem6.4}.
If the permutation $\sigma$ is even then the sets $Dg$ and  $Da$
are the same set, and the result is obvious.
Hence assume that $\sigma$ is an odd permutation of order $2^k$
and that $g$ has order $2^{k+1}$.

First we consider the case when
$k = \ell$ and $\sigma$ is a product of
$m$ disjoint $2^\ell$-cycles.  Since
changing an even number of signs in the representation of $g$ gives
and element that remains
in the coset $Dg$, we may assume that the first block in $g$
 is a $2^{\ell} \times 2^\ell$ matrix of the form
 \[C^{-} = \matfive{ 0 & 0& \cdots & 0 & -1\\
 1 & 0 & \cdots & & 0\\
 0 & 1 & 0 & \cdots & 0\\
\vdots & \ddots& \ddots &\ddots &\vdots\\
0 &\cdots  &0 &1 & 0},\]
the other blocks are the $2^\ell \times 2^\ell$ permutation matrix
 \[C^{+} = \matfive{ 0 & 0& \cdots & 0 & 1\\
 1 & 0 & \cdots & & 0\\
 0 & 1 & 0 & \cdots & 0\\
\vdots & \ddots& \ddots &\ddots &\vdots\\
0 &\cdots  &0 &1 & 0},\]
and that $a$ is the matrix with $m$ blocks of the permutation
matrix $C^{+}$.  Every element in $Dg$ is of order $2^{\ell +1}$,
since any even number
of changes of signs to $g$ will leave an odd number of sign changes
in at least one block.
Since $a$ has order $2^\ell$, then $Da$ has fewer elements of order
$2^{\ell + 1}$.

Next consider the remaining case: where $\sigma$ is an odd
permutation containing
at least one $2^k$-cycle and at least one $2^{t}$-cycle for
$0 \leq t < k$.
Without loss of generality we may assume that $g$ is
represented by a matrix with
top block a $2^k \times 2^k$ matrix of the form $C^-$, and
the final block is of
strictly smaller size.  Hence $a$ is represented by a matrix with
a top block a $2^k \times 2^k$ matrix of the form $C^+$, and
with smaller size matrices below.
If $g$ has $s$ blocks of order $2^k$ let $d=(d_1, \cdots, d_s, d_{s+1})$
represent a diagonal matrix where the $d_1, \cdots, d_s$ are each $2^k$-tuples
on the diagonal of $d \in D$.
Elements of $G$ of order $2^{k+1}$ in $Dg$ are of the form $da$, where
$d \in D$ has at least one of the $d_1, \cdots, d_s$ with an odd number of $-1$'s.
To establish  the result, we produce an injective map from the elements
of $G'$ of order $2^{k +1} $ in $Da$ to elements of $G$ of order
$2^{k +1}$ in $Dg$.
If at least one of $d_2, \cdots, d_s$ has an odd number of $-1$'s, then
associate $dg$ to $da$.  If none of $d_2, \cdots, d_s$ has
an odd number of $-1$'s
then $d_1$ has an odd number of $-1$'s so we associate
$d' dg$ to $da$, where $d'$ is the diagonal matrices with exactly two
$-1$'s (one in the first entry, and one in the last entry).
\end{proof}

Now we use the lemmas to show that the groups $G$ and $G'$ are not
isomorphic when $n$ is even.   We will show that $G$ has more
elements than $G'$ of order $2^{t}$
 for some $t$.  The largest possible order of elements that are a power of $2$  both in $G$ and
in $G'$ is $2^{r +1}$, where $2^r$ is the largest power of $2$
that is $\leq n$. This happens only when $Dg$ has order
$2^r$ in $G/D\cong S_n$.  By Lemma \ref{xxlem6.7}
the number of such elements in $G$ is $\geq$ the number of such
elements in $G'$.  If $G$ has more such elements than $G'$, we are done, so
assume that  both $G$ and $G'$ have the same number of elements
of largest possible order $2^{r +1}$; hence these cosets also contain
the same number of elements of order $2^r$.  Now consider all
elements of order
$2^r$ in $G$ and $G'$; by Lemma \ref{xxlem6.6} these elements
of $G$ are in the cosets $Dg$, where either $Dg$
has order $2^{r -1}$ or $2^r$ in $G/D$, but
we are assuming the number of elements of order $2^r$ in $G/D$
and $G'/D$ is the same, so
we need consider only elements of $G/D$ of order $2^{r-1}$.
Again by Lemma \ref{xxlem6.7} the number
of elements of $G$ of order $2^r$ is $\geq $ the number of
elements from $G'$
of order $2^{r}$.  Again if these numbers are different, we are done, and
hence without loss of generality we may assume that $G$ and $G'$ have
the same number of elements of order $2^t$ for $t \geq \ell +2$,
where $n = 2^\ell m$ for $m$ odd.  Now
consider elements of order $2^{\ell +1}$.
These elements arise only from $Dg$ with $Dg$ of order $2^\ell$ or
$2^{\ell +1}$.  Again we may assume that number of elements of
order $2^{\ell+1}$
arising from cosets in $G/D$ of order $2^{\ell +1}$ is the same for
$G$ and $G'$ because the number of elements of order $2^{\ell +2}$
in these cosets is
the same.  Finally consider cosets of order $2^\ell$.
For each of these cosets
the number of elements of order $2^{\ell +1}$ in $G$ is
$\geq$ the number of elements
of order $2^{\ell +1}$ in $G'$.  But we have the element of $G$ that has $m$
blocks of size $2^\ell$, for $m$ odd, so this coset is represented by
an element $g$ that has top block a $2^\ell \times 2^\ell$ block of the form
$C^-$ and all other blocks $2^\ell \times 2^\ell$ blocks $C^+$.  All elements
in this coset have order $2^{\ell +1}$, while the corresponding coset
of $D$ in $G'$ has the element with all $2^\ell \times 2^\ell$
blocks of the form $C^+$, so
this coset of $g'$ has at least one element has order $2^\ell$.
Hence $G$ and $G'$ have different
numbers of elements of order $2^{\ell +1}$, and these groups
are not isomorphic, establishing the claim that all $M(n,1,2)$ for
$n$ even are not isomorphic to a classical reflection group.
\end{example}

Forgetting about the underlying quantum polynomial ring $A$,
we can define a mystic reflection of a $k$-vector space.

\begin{definition}
\label{xxdefn6.8}
Let $V$ be a finite dimensional vector space over $k$.
A linear map $g: V\to V$ is called a {\it mystic reflection}
if there is a basis of $V$, say $\{y_1,\cdots,y_n\}$,
such that $g(y_j)=y_j$ for all $j\leq n-2$ and $g(y_{n-1})
=i\; y_{n-1}$ and $g(y_n)=-i\; y_n$.
\end{definition}

Theorem \ref{xxthm0.1} says that every reflection group
of a skew polynomial ring is generated by reflections
and mystic reflections of the vector space $A_1$.
Further such a reflection group is a product of
classical reflection groups of $A_1$ and copies of
$M(n,\alpha,\beta)$'s. We finish the paper by the following
question.

\begin{question}
\label{xxque6.9}
Let $G$ be a finite subgroup of $GL(V)$. If $G$ is
generated by reflections in the classical sense and
mystic reflections in the sense of Definition \ref{xxdefn6.8},
then is $G$ isomorphic to  a product of
classical reflection groups of $V$ and
copies of $M(n,\alpha,\beta)$'s?
\end{question}

\section*{Acknowledgments}
J.J. Zhang is supported by the National Science Foundation of USA and
the Royalty Research Fund of the University of Washington.



\begin{thebibliography}{10}



\bibitem[BG]{BG}
K. A. Brown and K. R. Goodearl,
Lectures on Algebraic Quantum Groups, Birkh{\"a}user, 2002.


\bibitem[CLO]{CLO} D. Cox, J. Little, D. O'Shea,
Ideals, Varieties, and Algorithms: An Introduction
to Computational Algebraic Geometry and Commutative
Algebra (Second Edition), Springer, New York, 1997

\bibitem[JiZ]{JiZ}
N. Jing and J.J. Zhang,
On the trace of graded automorphisms,
J. Algebra {\bf 189} (1997), no. 2, 353--376.

\bibitem[JoZ]{JoZ}
P. J{\o}rgensen and J.J. Zhang,
Gourmet's guide to Gorensteinness,
Adv. Math. {\bf 151} (2000), no. 2, 313--345.

\bibitem[K]{K} R. Kane, Reflections Groups and Invariant
Theorem, CMS Books in Mathematics {\bf 5}, Springer, New York, 2001.

\bibitem[KKZ1]{KKZ1}
E. Kirkman, J. Kuzmanovich and J.J. Zhang,
Rigidity of graded regular algebras, Trans. Amer. Math. Soc.
(to appear)

\bibitem[KKZ2]{KKZ2}
E. Kirkman, J. Kuzmanovich and J.J. Zhang,
Hopf algebra (co)-actions on Artin-Schelter regular algebras,
in preparation.


\bibitem[N]{N}
M.D. Neusel,
Invariant Theory,
Student Mathematical Library, Vol {\bf 38},
American Mathematical Society, Providence, R.I., 2007.

\bibitem[ShT]{ShT}
G.C. Shephard and J.A. Todd,
Finite unitary reflection groups,
Canadian J. Math. {\bf 6}, (1954). 274--304.


\bibitem[YZ]{YZ}
A. Yekutieli and J.J. Zhang,
Homological transcendence degree,
Proc. London Math. Soc. (3) {\bf 93} (2006),
no. 1, 105--137.

\bibitem[Z]{Z}
J.J. Zhang,
Twisted graded algebras and equivalences of graded categories,
Proc. London Math. Soc. (3) {\bf 72} (1996), no. 2, 281--311.


\end{thebibliography}
\end{document}